\documentclass[a4paper,11pt,oneside,reqno]{amsart}

\usepackage{amsmath}
\usepackage{amsthm}
\usepackage{amssymb}
\usepackage{graphicx}

\theoremstyle{plain}
\newtheorem{teore}{Theorem}[section]
\newtheorem{defin}[teore]{Definition}
\newtheorem{lem}[teore]{Lemma}
\newtheorem{coro}[teore]{Corollary}
\newtheorem{propo}[teore]{Proposition}
\newtheorem*{claim}{Claim}
\newtheorem*{claim*}{Claim}
\theoremstyle{remark}
\newtheorem{ejemplo}[teore]{{\sc Example}}
\newtheorem{notas}[teore]{{\sc Remark}}

\newcommand{\prop}{\begin{propo}}
\newcommand{\fprop}{\end{propo}}
\newcommand{\cor}{\begin{coro}}
\newcommand{\fcor}{\end{coro}}

\newcommand{\defi}{\begin{defin}\rm}
\newcommand{\fdefi}{\end{defin}}
\newcommand{\eje}{\begin{ejemplo}}
\newcommand{\feje}{\end{ejemplo}}
\newcommand{\lema}{\begin{lem}}
\newcommand{\flema}{\end{lem}}
\newcommand{\teor}{\begin{teore}}
\newcommand{\fteor}{\end{teore}}
\newcommand{\nota}{\begin{notas}\rm}
\newcommand{\fnota}{ \end{notas}}
\newcommand{\clam}{\begin{claim}}
\newcommand{\fclam}{\end{claim}}
\newcommand{\clams}{\begin{claim*}}
\newcommand{\fclams}{\end{claim*}}

\newcommand{\lclam}{\begin{lclaim}}
\newcommand{\flclam}{\end{lclaim}}
\newcommand{\prucl}{\prue[Proof of Claim:]}
\newcommand{\fprucl}{\fprue}
\newcommand{\ben}{\begin{enumerate}}
\newcommand{\een}{\end{enumerate}}
\newcommand{\bit}{\begin{itemize}}
\newcommand{\eit}{\end{itemize}}
\newcommand{\mc}[1]{\mathcal{#1}}

\newcommand{\casos}{\begin{itemize}}
\newcommand{\fcasos}{\end{itemize}\setcounter{cs}{1}}

\newcommand{\fin}{\textsc{FIN}}
\newcommand{\pe}{\preceq}

\newcommand{\mU}{\mathcal{U}}

\renewcommand{\iff}{if and only if }

\DeclareMathOperator{\Rel}{R}

\newcommand{\conj}[2]{ \{ {#1}\,:\,{#2} \} }

\newcommand{\buit}{\emptyset}

\newcommand{\ga}{\gamma}
\newcommand{\A}{\mathbb{A}}
\newcommand{\B}{\mathbb{B}}

\DeclareMathOperator{\rel}{\thicksim}
\DeclareMathOperator{\nrel}{\not\thicksim}
\newcommand{\rang}{\text{Range }}
\newcommand{\Ga}{\Gamma}

\newcommand{\al}{\alpha}
\newcommand{\be}{\beta}
\newcommand{\de}{\delta}

\newcommand{\linspan}{\text{LinSpan }}
\newcommand{\La}{\Lambda}

\newcommand{\vep}{\varepsilon}

\newcommand{\R}{{\mathbb R}}

\newcommand{\N}{{\mathbb N}}

\newcommand{\med}{\mathrm{ mid}}
\newcommand{\mini}{\mathrm{ min}}
\newcommand{\maxi}{\mathrm{ max}}

\newcommand{\rest}{\negmedspace\negmedspace\upharpoonright\negthickspace}

\newcommand{\dom}{\text{dom }}

\newcommand{\supp}{\text{supp }}

\newcommand{\con}{\subseteq}

\newcommand{\prue}{\begin{proof}}
\newcommand{\fprue}{\end{proof}}

\begin{document}
\title{Canonical equivalence relations on nets of $PS_{c_0}$}
\author{J. Lopez-Abad}
\address{ Equipe de Logique Math\'{e}matique \\
Universit\'{e} Paris 7- Denis Diderot\\
C.N.R.S. -UMR 7056\\
2, Place Jussieu- Case 7012\\
75251 Paris Cedex 05 \\
France}\email{abad@logique.jussieu.fr} 
\thanks{This Research was supported through a
European Community Marie Curie Fellowship}%
\begin{abstract}
We give a list  of canonical  equivalence relations on discrete nets of the positive unit
sphere of $c_0$. This  generalizes  results of W. T. Gowers \cite{gow2} and A. D. Taylor
\cite{tay}.
\end{abstract} \maketitle

\section{Introduction}
Let $\fin$ be the family of nonempty finite sets of positive integers. A block sequence is an
infinite sequence $(x_n)_n$ of elements of $\fin$ such that  for every $n$ one has $\max x_n<\min
x_{n+1}$ (usually written as $x_n<x_{n+1}$). The \emph{combinatorial subspace} $\langle (x_n)_n
\rangle$ given by   $(x_n)_n$  is  the set of finite unions $x_{n_0}\cup \cdots \cup x_{n_m}$.
Using this terminology,  the Hindman's pigeonhole principle   \cite{hind} of $\fin$ states that
every finite coloring of $\fin$ is constant in some combinatorial subspace, or, equivalently, every
equivalent relation on $\fin$ with finitely many classes   has a restriction to
 some   combinatorial subspace with only one class. It is easy to see, for example by
considering the equivalent relation defined by $s\rel t$ iff $\min s=\min t$, that this is no
longer the case for equivalence relations with an arbitrary number of classes. Nevertheless, it is
still possible to classify them, much in the spirit of the original motivation of F. P. Ramsey
\cite{ram} for discovering his famous Theorem. A result of  Taylor \cite{tay} states that every
equivalence relation on $\fin$ can be reduced, by restriction to a combinatorial subspace, to one
of the following five canonical relations:
\begin{equation*}
\min,\max,(\min,\max), =, \fin^2,
\end{equation*}
naturally defined by   $s \min t$ iff  the minimum of $s$ is equal to the minimum of $t$, $s \max
t$ iff the maximum of $s$ is equal to the maximum of $t$, $s (\min,\max)t$ iff both minimum and
maximum are the same.

Following some geometric ideas   exposed in Section 2, one can   generalize $\fin$ as follows:
Given a positive integer $k$, let $\fin_k$ be the set of mappings $x:\N\to \{0,1,\dots,k\}$, called
\emph{$k$-vectors}, whose support $\supp x=\conj{n}{x(n)\neq 0}$ is finite and with $k$ in their
range.  One can naturally extend the union operation on $\fin$ to the join operation $\vee$ on
$\fin_k$   by $(x\vee y )(n)=\max\{x(n),y(n)\}$.  Let $T:\fin_k\to \fin_{k-1}$ be the mapping
defined by $T(x)(n)=\max\{x(n)-1,0\}$. A $k$-block sequence $(x_n)_n$ is an infinite sequence of
members of $\fin_k$ such that $\max\supp x_n<\min\supp x_{n+1}$ for every $n$.  The
$k$-\emph{combinatorial subspace} $\langle (x_n)_n \rangle$ defined by a $k$-block sequence $(x_n)$
is the set of combinations of the form $T^{i_0}x_{n_0} \vee\cdots \vee T^{i_m}x_{n_m}$ with the
condition that $i_j=0$ for some $j$, and where  $T^{i}x$ is defined by $T^i x (n)=\max\{x(n)-i,0\}$
for $i>0$ and  $T^0=Id$.  Gowers   has proved in \cite{gow2} that $\fin_k$ possesses the exact
analogue of the pigeonhole principle of $\fin$: Every equivalence relation on $\fin_k$ with
finitely many classes has a restriction to
 some   combinatorial subspace with only one class.
The aim of this paper is to characterize  equivalence relations on $\fin_k$ with arbitrary number
of classes. More precisely, we are going to give a non redundant finite list $\mathcal{T}_k$ of
equivalence relations such that any other equivalence relation on $\fin_k$ can be reduced, modulo
restriction to some $k$-combinatorial subspace, to one in the list $\mathcal{T}_k$.

Indeed,   the elements  of  $\mc T_k$ are determined by   characteristics of a typical $k$-vector.
Easy examples of these are the minimum and maximum of a finite set, that determine the Taylor's
list for $\fin$. Generalizing this, given an  integer $i$ with $1\le i \le k$  let $\min_i s$ be
  the least integer $n$ such that $s(n)=i$. Another more complex example is the following. Given two integers $i$ and $l$ such that
$1\le l\le i-1\le k$,  let us assign to a given vector $s$ of $\fin_k$ the set of integers $n$ such
that $\min_{i-1}s\le n\le \min_{i}s$ and $s(n)=l$.  We illustrate this with the following picture.
\begin{center}
\begin{figure}[h]
\includegraphics[scale=1]{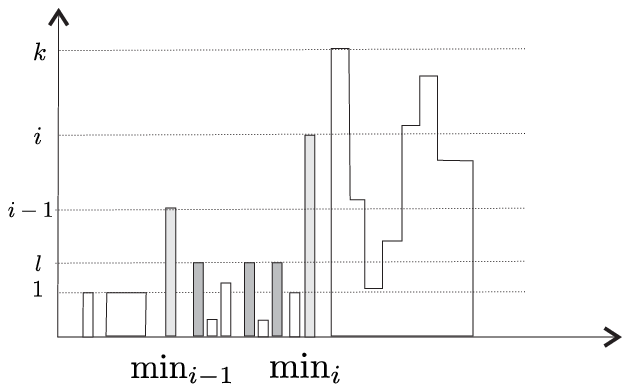}
\caption{ an example of invariant}
\end{figure}
\end{center}
So, our first task will be to guess all the natural characteristics of a $k$-vector. Although these
characteristics  are not well defined for an arbitrary $k$-vector, we will show that every
$k$-block sequence will have a $k$-\emph{block} subsequence, called here a \emph{system of
staircases} for which all the vectors have all natural characteristics well defined. The precise
definitions are given in Section 3.

In order to show that every equivalence relation is, when restricted to some $k$-combinatorial
subspace, in  $\mc T_k$ we follow the ideas  of Taylor's proof
  \cite{tay}.  Let us explain this.  Given an equivalence relation $\rel$ on $\fin$ one defines the coloring $c:[\fin]^{[3]}\to
\{0,1\}^{4}$  by
\begin{equation}\label{taylora}
a=(a_0,a_1,a_2)\mapsto \left\{\begin{array}{ll} c(a)(0)=1 &\text{iff } a_0\rel
a_1 \\
c(a)(1)=1 &\text{iff } a_0\cup a_1\rel a_0 \\c(a)(2)=1 &\text{iff } a_0\cup a_1\rel a_1
\\c(a)(3)=1 &\text{iff } a_0\cup a_1\cup a_2\rel
a_0\cup a_2,
\end{array} \right.
\end{equation}  where $[\fin]^{[3]}$ is the set of 3-sequences of finite sets $(a_0,a_1,a_2)$ such that
$a_0<a_1<a_2$. Since $[\fin]^{[3]}$ has a pigeonhole principle (this is a simple extension of
Gowers' result), one  can find a block sequence $X=(x_n)_n$ such that $c$ is constant on
$[X]^{[3]}$ with value $s_0\in \{0,1\}^4$.  An analysis  of the value  $s_0$ identifies the
restriction of the equivalence relation $\sim$ to $X$ as one of the five relations
$\min,\max,(\min,\max),=,\fin^2$.  Let us  re-write  the coloring $c$  in a way that will be easy
to generalize to $\fin_k$. Fix an alphabet of countably many variables $\{x_n\}_n$. An
$\rel$-equation $e$ is a pair $((x_{i_0},\dots,  x_{i_l}),(x_{j_0},\dots, x_{j_m}))$,   written as
$ x_{i_0}\cup \cdots \cup x_{i_l}\rel x_{j_0}\cup \cdots\cup x_{j_m}$, such that $0=i_0<\cdots
<i_l$, $j_0<\cdots <j_m$. We say that  equation $e$ is true in $X$ iff for every sequence
$a_{0}<\cdots<a_{\max\{{i_l,j_m}\}}$ in $X$ the corresponding substitutions  $a_{i_0}\cup \cdots
\cup a_{i_l}$ and $ a_{j_0}\cup \cdots\cup a_{j_m}$ are $\rel$-related; we say that the equation
$e$ is false in $X$ iff for every sequence $a_{0}<\cdots<a_{\max\{{i_l,j_m}\}}$ in $X$ one has that
$a_{i_0}\cup \cdots \cup a_{i_l}\nrel a_{j_0}\cup \cdots\cup a_{j_m}$. The equation $e$ is
\emph{decided} in $X$ if it is either true or false in $X$. Using this terminology,  one can
re-state the fact that the coloring $c$ is constant on $X$  by saying  that the equations $x_0\rel
x_1$, $x_0\cup x_1\rel x_0$, $x_0\cup x_1\rel x_1$ and $x_0\cup x_1\cup x_2\rel x_0\cup x_2$ are
all decided in $X$. Taylor's proved that these four equations determine the equivalence relation
$\rel$. For an arbitrary integer $k$, the list of equations to be considered is, obviously, longer.
For example, for $k=2$ the equations
\begin{equation*}
x_0\cup x_1\cup Tx_2\sim x_1\cup Tx_2\text{ and }x_0+Tx_1\cup x_2\sim x_0\cup x_2,
\end{equation*}
need to be  considered. So, the next goal, after one has identified the list $\mc T_k$, is to find
a set $L$ of $\sim$-equations  characterizing a given equivalence relation $\sim$ on $\fin_k$. The
first candidate for $L$ is the set of all equations. It turns out  that the lists $\mc T_k$
consists on all the equivalence relations for which every equation is always true or always false,
independently of the $k$-block sequence considered.  So  it does not seem   reasonable to try to
find directly a $k$-block sequence deciding all equations. Instead, we   first find a  smaller list
of equations
 decided in some $k$-block sequence, but at the same time large enough  to   use the inductive hypothesis to provide
a richer list of equations, determining our given equivalence relation as one of the list $\mc
T_k$.

It is worth to point out that we give an explicit description of $\mc T_k$ in a way that it is
possible to describe the number  $t_k$ of equivalence relations in $\mathcal{T}_k$  using standard
arithmetic functions, as for example the incomplete $\Gamma$ function:
\begin{equation*}
t_k=|\mathcal{T}_k|=e^2 \left[ k\left[\Gamma(k, 1) - \Gamma(k+1, 1) \right]^2 +
        \Gamma(k+1,1)^2\right].\end{equation*}
Since  $\fin_k$   is isomorphic to a net of  the positive sphere of $c_0$, our result implies the
immediate analogue for those nets. For example, given an equivalence relation $R$ on $PS_{c_0}$ and
given some $\de>0$ there is an infinite dimensional block subspace $X$ of $c_0$ and some
equivalence relation $R'$ in our finite list such that every  $R'$-class in $X$ is included in the
$\de$-fattening of some $R$-class.

This paper is organized as follows. In Section 2 we introduce $\fin_k$ as a natural copy of a net
of the positive sphere of $c_0$,  extending some standard concepts coming from Banach space theory
to $\fin_k$. We also state the W. T. Gowers Pigeonhole principle of  $\fin_k$. The notion of
equation is introduced in Section 3, together with the natural characteristics of a vector of
$\fin_k$. We describe the vectors for which these invariants are well defined, and we show that
they appear ``everywhere". We also define  the family $\mathcal{T}_k$. In Section 4  our main
theorem is proved, and in Section 5 we give an explicit formula to compute the cardinality of
$\mathcal{T}_k$. Sections 6 and 7 deal with the finite version of our main result, and with some
consequences for equivalence relations on the positive sphere of $c_0$.

\section{First Definitions and Results}

Recall that $c_0=c_0(\R)$ is the Banach space of sequences of real numbers  converging to 0, with
the $\sup$-norm defined for a vector $\vec{x}=(x_n)_n$ of $c_0$ by $\|\vec{x}\|=\sup_n|x_n|$. Let
$(e_n)_n$ be its natural Schauder basis, i.e., $e_n(m)=\de_{n, m}$. The support of a vector
$\vec{x}=(x_n)_n$, is defined as $\supp \vec{x}=\conj{n}{x_n\neq 0}$ and let $c_{00}$ be the linear
subspace of $c_0$ consisting of the vectors $\vec{x}=(x_n)_n$ with finite support, i.e., only
finitely many of the coordinates of $\vec{x}$ are not zero. Given two vectors $\vec{x}$ and
$\vec{y}$ of $c_{00}$ we write $\vec{x}<\vec{y}$ to denote that $\max\supp \vec{x}<\min \supp
\vec{y} $.

Let $PS_{c_0}$ be the set of norm one positive vectors of $c_0$, i.e., the set of all vectors
$\vec{x}=(x_n)_n$ such that $\|\vec{x}\|=1$, and such that $x_n\ge 0$, for every $n$, and let
$PB_{c_0}$ be the set of positive vectors of the unit ball of $c_0$. Observe that $PB_{c_0}$ is a
lattice with respect to $(x_n)_n\vee (y_n)_n=(\max\{x_n,y_n\})_n$ and $(x_n)_n\wedge
(y_n)_n=(\min\{x_n,y_n\})_n$, with $0=(0)_n$, and $1=(1)_n$.  Notice also that $PS_{c_0}$ is closed
under the operation $\vee$, and that $x\vee y=x+y$ if $x$ and $y$ have disjoint support. In
general, given two subsets $N\con A$ of $c_0$ and a positive number $\de$ we say that $N$ is a
\emph{$\de$-net} of $A$ iff for every $\vec{a}\in A$ there is some $\vec{x}\in N$ such that
$\|\vec{a}-\vec{x}\|\le \de$.

For a given $\de$ with $0<\de <1$,  let $k$ be the least integer such that $1/(1+\de)^{k-1}\le
\de$, and let $\vep=1/(1+\de)$.
 Let
\begin{align*}
\mathcal{N}_{\de}=&\conj{x \in PB_{c_{00}}}{   x(i)\in \{1,\vep,\vep^2,\dots,\vep^{k-1},0 \} } \\
\mathcal{M}_{\de}=& \conj{x \in PS_{c_{00}}}{   x(i)\in \{1,\vep,\vep^2,\dots,\vep^{k-1},0 \} }.
\end{align*}
Since $\vep^i-\vep^{i+1}=\vep^i(1-\vep)=\vep^i(\de/(1+\de))< \de$ and $\vep^{k-1}\le \de$, it
follows that $\mathcal{N}_{\de}$ and $\mathcal{M}_\de$ are  $\de$-nets of $PB_{c_0}$  and
 of $PS_{c_0}$, respectively. The set  $\mathcal{N}_\de$ is a sub-lattice of $PB_{c_0}$ with respect to
$\vee$ and $\wedge$,  and  it is closed under scalar multiplication by $\vep$, identifying
$\vep^l=0$ for $l\ge k$ (which means that we identify the coordinates less than $\vep^k$ with 0).
Also, for two $\vec{x},\vec{y}\in \mathcal{M}_\de$, we have that $\vec{x}\vee\vep^i \vec{y}, \vep^i
\vec{x}\vee \vec{y} \in \mathcal{M}_\de$, for every $0\le i \le k-1$. Finally, note that
$\mathcal{N}_\de=\bigcup_{i=0}^{k-1} \vep^i \mathcal{M}_\de$, a disjoint union.

We define the mapping $\Theta=\Theta_\de:\mathcal{N}_\de\to \{0,1,\dots,k\}^\N$ by
\begin{equation*}
\Theta((x_m)_m)(n)=\left\{\begin{array}{ll} k-\log_\vep (x_n) & \text{if } x_n\neq 0 \\
0 & \text{if }x_n=0.
\end{array}\right.
\end{equation*}
We can now give an equivalent definition of $\fin_k$  using the mapping $\Theta$, and therefore giving a
geometrical interpretation of it.
\defi Fix, for a given integer $k$, a positive real number $\de=\de(k)$ such that
  $1/(1+\de)^{k-1}=\de$. Let
$\fin_k=\Theta(\mathcal{M}_\de)$, i.e.,   the set of functions $s:\N\to \{0,1,\dots,k\}$ eventually
0, and with $k$ in the range.  The elements of $\fin_k$ are called $k$\emph{-vectors}.

Observe that $\Theta"\mathcal{N}_\de=\bigcup_{i=0}^{k-1}\Theta" \vep^i \mathcal{M}_\de$, and that
$\Theta" \vep^i \mathcal{M}_\de$ is the set of all functions $s:\N\to \{0,1,\dots,k\}$ eventually
0, and with $k-i$ in the range. So,  $\Theta" \vep^i \mathcal{M}_\de=\fin_{k-i}$. Hence,
$\Theta"\mathcal{N}_\de=\bigcup_{i=1}^{k}\fin_i=: \fin_{\le k}$, whose members are called
\emph{$(\le \!k)$-vectors.}
\fdefi
 We can transfer the algebraic structure of $\mathcal{N}_k\con
c_{00}$ to $\fin_{\le k}$ via $\Theta$. In particular,  for $s,t\in \fin_{\le k}$, let the \emph{support of
$s$} be $\supp s=\conj{n}{s(n)\neq 0}$; we write $s<t$ to denote that $\max \supp s <\min\supp t$, define
$s\vee t$ and $s\wedge t$ by
\begin{equation*}
(s\vee t)(n)=\max \{s(n),t(n)\} \text{ and }(s\wedge t)(n)=\min \{s(n),t(n)\},
\end{equation*}
and let $T$ be the transfer of the multiplication by $\vep$,   i.e., for a $(\le\! k)$-vector $s$,
let
\begin{equation*}
T(s)=\Theta(\vep\Theta^{-1} (s))=(s-1)\vee 0.
\end{equation*}
Let $S:\fin_{k-1}\to \fin_k$ be an inverse map for $T$, defined for a $(k-1)$-vector $a$ by
\begin{equation*}
S(a)(n)=\left\{
\begin{array}{ll}
a(n)+1 & \text{if } n\in \supp a\\
0 & \text{if not.}
\end{array} \right.
\end{equation*}
It turns out that $\fin_{\le k}$ is a lattice with operations $\vee$ and $\wedge$, and it is closed
under $T$. We will use the order $\le_L$ to denote the lattice-order of $\fin_{\le k}$, i.e., for
$s,t\in \fin_{\le k}$, we write $s\le_L t$ iff $s\wedge t=s$. Note  that $\fin_i \vee \fin_j
=\fin_{\max_{\{i,j\}}}$ and $\fin_i \wedge\fin_j =\fin_{\min_{\{i,j\}}}$.  We will use $s+t$ for
$s\vee t$  whenever $s<t$.

We now pass to introduce some combinatorial notions. A sequence of $k$-vectors $(s_n)$ is called a
 \emph{finite $k$-block sequence} if  $(s_n)$ is finite and if  $s_n<s_{n+1}$ for every $n$; if such sequence is infinite, then
 we call it a (infinite) \emph{$k$-block sequence}. We write  $\fin_k^{[\infty]}$, $\fin_k^{[n]}$ and $\fin_k^{[<\infty]}$ to denote
 respectively
the set of  $k$-block sequences, finite  $k$-block sequences of length $n$, and the set of finite $k$-block
sequences.

The  \emph{$k$-combinatorial subspace}   $\langle \alpha \rangle$ defined  by  a finite or infinite $k$-block
sequence  $\alpha=(s_n)_{n}$  is  the set of all $k$-vectors of $\al$ defined by
$$\langle \al \rangle=\Theta((\linspan\Theta^{-1}\{s_n\}_n )\cap \mathcal{N}_\de),$$
where $\linspan A $ denotes the linear span of a given subset $A$ of $c_0$. Using this one has that
$\fin_k=\langle (\Theta e_n)_n \rangle$. Similarly, we define for  a given  integer $i\le k$ the
set $\langle \al \rangle_i$  of $i$-vectors of $\al$. A main property of the $k$-block sequences
$(a_n)_n$ is that $e_n\mapsto a_n$ naturally extends to a lattice isomorphism between $\fin_k$ and
$\langle (a_n)_n\rangle$ that preserves the operation $T$.

 For $M\le N\le \infty$, and $\al=(s_n)_{n<N}$ let $[\al]^{[M]}$ be the
set of \emph{$k$-block subsequences} of $\al$, defined as $[\al]^{[M]}=\conj{(s_n)_{n<M}\in
\fin_k^{[M]}}{s_n\in \langle \al \rangle\, (0\le n < M)}$. Without loss of generality we will identify
$[\al]^{[1]}$ with $\langle \al \rangle$.

Given two finite block sequences $\al$ and $\be$, and two infinite ones $A$ and $B$, we define
$\al\pe \be$ \iff $\al\in [ \be ]^{[|\al|]}$, $\al\pe A$ \iff $\al\in [A]^{[|\al|]} $, and  $B\pe
A$ \iff $B\in [A]^{[\infty]}$. Notice that all these definitions come from the notion of subspace.
For example, $A\in \langle B\rangle$ \iff the space generated by $\Theta^{-1} A $ is a subspace of
the space generated by $\Theta^{-1} B $.

For a $k$-block sequence  $A=(a_i)_i$ and $a\in\langle A\rangle$, since $\langle A
\rangle=\Theta(\linspan\Theta^{-1}\{a_i\}_i  \cap \mathcal{N}_\de)$, we have that  $\Theta^{-1}a\in
\Theta^{-1}\{a_i\}_i \cap \mathcal{N}_\de $. Therefore, $\Theta^{-1}a=\sum_{i=0}^m \vep^{d_i}
\Theta^{-1}a_i$, for some $m$, and with possibly some $d_i=0$. This implies that
$a=\Theta(\sum_{i=0}^m \vep^{d_i} \Theta^{-1}a_i)=\sum_{i=0}^m \Theta(\vep^{d_i}
\Theta^{-1}a_i)=\sum_{i=0}^m T^{d_i}a_i$.

Finally, an infinite sequence $(A_r)_{r\in \N}$ of infinite $k$-block sequences $A_r=(a_n^r)_n$ is
called a \emph{fusion sequence  of $A\in \fin_k^{[\infty]}$} if for all $r\in \N$:

\noindent (a) $A_{r+1}\pe A_r\pe A$,

\noindent (b) $a^r_0<a^{r+1}_0$.

The infinite $k$-block sequence  $A_{\infty}=(a_0^r)_{r\in \N}$ is called the  \emph{fusion
$k$-block sequence}  of the sequence $(M_r)_{r\in \N}.$

\defi
Given a $k$-block sequence $A=(a_n)_n$, let $C_A:\langle A\rangle\to FIN_k$ be  the mapping
satisfying
\begin{equation}\label{mnjuew}
a=\sum_{n=0}^\infty T^{k-C_A(a)(n)}a_{n},
\end{equation} for every $k$-block vector $a$ of $A$.  Since $\Theta^{-1}a=\sum_{n\ge 0}\vep^{k-C_A(a)(n)} \Theta^{-1}a_n$,
for every $a$, the mapping  $C_A$ is well defined. We call the sum in (\ref{mnjuew}) the
\emph{canonical decomposition of $a$ in $A$}. Notice that $C_A(a)\in \fin_k$ for every $a$.

For two $(\le\! k)$-vectors $s$ and $t$,

\noindent (a) we write $s\sqsubseteq t$ when $t\rest \supp s =s$, i.e., if $t$ restricted to the
support of $s$ is equal to $s$, and

\noindent (b) we write $s\perp t$ when there is no $u\in \fin_{\le k}$ such that $u\sqsubseteq
s,t$, i.e., if $s(n)\neq t(n)$ for every $n\in \dom s\cap \dom t$.

Using this, if $s=\sum_{n=0}^\infty T^{k-l_n}a_{n}$, then $T^{k-l_n}a_{n}\sqsubseteq a$, for every
$n$, while $T^{k-l_n}a_{n}\perp T^{k-l_{n'}}a_{n'}$ for every $n\neq n'$. It follows that:
\fdefi
\prop\label{basis}
Fix $A=(a_n)_n$,  $a\in \langle A \rangle$ and an integer $n$. If there are some  $r\le k$  and $m$
such that $T^{k-r}a_n (m)=a(m)\neq 0$, then necessarily $C_A(a)(n)=r$ (i.e., $T^{k-r}a_n
\sqsubseteq a$). \qed
\fprop

The following is Gowers' pigeonhole principle for $\fin_k$.
\teor\cite{gow2}\label{gow}
If $\fin_k$ is partitioned into finitely many pieces, then there is $A\in
\fin_k^{[\infty]}$ such that $\langle A \rangle$ is in only one of the pieces.
\fteor

This naturally extends to higher dimensions.
\lema \label{gowgen} \cite{tod} Suppose that $f:\fin_k^{[n]}\to \{0,\dots,l-1\}$. Then there is a block
sequence $X$ such that  $f$ is constant on $[X]^{[n]}$.\flema
\prue
The proof is done by induction on $n$. Suppose it is true for $n-1$. We can find, by a repeated use
of Theorem \ref{gow}, a fusion sequence $(X_r)_r$, $X_r=(x_i^r)_i$, such that for every $r$ and
every $(b_0,\dots,b_{n-2})\in [(x_i^i)_{i<r}]^{[n-1]}$ the coloring $f$ is constant on the set
 $\conj{(b_0,\dots,b_{n-2},x)}{x\in X_r}$ with value $\vep((b_0,\dots,b_{n-2}),r)$. By construction one has that  $X_r \pe
 X_{s}$ if $r\le s$. So it follows that
 $\vep((b_0,\dots.,b_{n-2}),
r)=\vep((b_0,\dots,b_{n-2}), s)$ for every $(b_0,\dots,b_{n-2})\in [\theta_r ]^{[n-1]}$ and every
$r<s$. This allows us to define  $\vep:[X_\infty]^{[n-1]}\to \{0,1,\dots,l-1\}$ by
$\vep(b_0,\dots,b_{n-2})=\vep((b_0,\dots,b_{n-2}),r)$, for some (any) integer $r$, where
$X_\infty=(x_i^i)_i$ is the fusion $k$-block sequence of $(X_r)_r$. This coloring $\vep$ can be
easily interpreted as a coloring of $\fin_{k}^{[n-1]}$, so by the inductive hypothesis there is
some $X\pe X_\infty$ such that $\vep$ is constant on $[X]^{[n-1]}$, and therefore $f$ is also
constant on $[X]^{[n]}$.
\fprue

\section{Equations, Staircases and Canonical equivalence Relations}

Roughly speaking, terms are natural mappings that assign $k$-vectors to finite block sequences of $k$-vectors
of a fixed length $n$, and which are defined from the operations $+$ and $T^i$ of $\fin_k$. For example, the
mapping that assigns to a block sequence $(a_1,a_2)$ of $k$-vectors the $k$-vector $a_1+Ta_2$ is a $k$-term
which can be understood as the mapping with two variables $x_1,x_2$ defined by $f(x_1,x_2)=x_1+T x_2$.

From   two fixed $k$-terms $f$ and $g$ of $n$ variables  and one equivalence relation $\sim$ on
$\fin_k$  we can define the natural coloring $c_{f,g}:[\fin_k]^{[n]}\to \{0,1\}$ via
$c_{f,g}(a_1,\dots,a_n)=1$ if and only if $f(a_1,\dots,a_n) \sim g(a_1,\dots,a_n)$. A $k$-equation
will be $f\sim g$. The pigeonhole principle in Lemma  \ref{gowgen} gives that for every equation
$f\sim g$ ($f$ and $g$ with $n$ variables) there is some infinite block sequence $A$  such that,
either for every $(a_1,\dots,a_n)$ in $[A]^{[n]}$, $f(a_1,\dots,a_n) \sim g(a_1,\dots,a_n)$, or for
all $(a_1,\dots,a_n)$ in $[A]^{[n]}$, $f(a_1,\dots,a_n)\not\sim g(a_1,\dots,a_n)$, i.e., in $A$ the
equation $f\sim g$ is either true or false. As we explained in the introduction, Taylor proves that
an equivalence relation $\sim$ on $\fin$ is determined by a list of 4 equations (precisely,
$x_0\sim x_1$, $x_0\sim x_0+x_1$, $x_1\sim x_0+x_1$ and $x_0+x_1+x_2\sim x_0+x_2$). This is going
to be also the case for arbitrary $k$, of course with a more complex list of equations.

\subsection{Terms and equations}
\defi
Let $\texttt{X}=\{x_n\}_{n\ge 1}$ be a countable infinite alphabet of  variables. Consider the
trivial map $\texttt{x}: \texttt{X} \to \N$, defined by $x_n  \mapsto \texttt{x}(x_n)=n$. A
\emph{free $k$-term} $\texttt{p}$ is a map of the form $s\circ \texttt{x}$ where $s$ is a
$k$-vector, i.e., it is a map $\texttt{p}:\texttt{X}\to \{0,\dots,k\}$ such that $\supp \texttt{p}$
is finite, and $k$ is in the range of $\texttt{p}$. A natural representation of $\texttt{p}$ is
$$\texttt{p}=\texttt{p}(x_0,\dots,x_l)=\sum_{i=0}^l T^{k-m_i}x_{i},$$
where $0\le m_i \le k$, and at least one $m_i=k$. For example $T^2 x_1+Tx_2+x_4$, and $x_1+x_5$ are both free
3-terms. Notice that, if $\texttt{p}$ is a free $k$-term, then $\texttt{p}\circ \texttt{x}^{-1}$ is a
$k$-vector. A \emph{free $(\le \! k)$-term} is $s\circ \mathbf{x}$, where $s$ is a $(\le \! k)$-vector. It
follows that the set of free $(\le\! k)$-terms is a lattice. For example
\begin{equation*}
\texttt{p}(x_0,\dots,x_n)\vee \texttt{q}(x_0,\dots,x_m)= (\texttt{p}\circ \texttt{x}^{-1} \vee
\texttt{q} \circ \texttt{x}^{-1})\circ \texttt{x}.
\end{equation*}
We also have  defined the operator $T$ for a $k$-term $p(x_0,\dots, x_n)$ by
$$T(p(x_0,\dots,x_n))=(T(p\circ \texttt{x}^{-1})\circ \texttt{x}).$$
\fdefi
For every  $(\le \! k)$-term $p(x_0,\dots,x_n)=\sum_{i=0}^n T^{k-m_i}x_{i}$ we consider the
following kind of  substitutions:

 \noindent (a)  Given a sequence of free $(\le \! k)$-terms
$t_0,\dots,t_n$, consider the substitution of each $x_i$ by $t_i$
\begin{equation*}
\texttt{p}(t_0,\dots,t_n)=\bigvee_{i=0}^n T^{k-m_i}t_{i}.
\end{equation*}
In the case that $\texttt{p}$ and $t_0,\dots,t_n$ are  free $k$-terms, then
$\texttt{p}(t_0,\dots,t_n)$ is also a free $k$-term.

\noindent (b) For a block sequence $(a_0,\dots,a_n)$ of $(\le \!k)$-vectors,  replace each $x_i$ by
$a_i$
\begin{equation*}
\texttt{p}(a_0,\dots,a_n)=\sum_{i=0}^n T^{k-m_i}a_{i}.
\end{equation*}
If $\texttt{p}$  is a free $k$-term, and $a_0,\dots,a_n$ are $k$-vectors, then the result of the
substitution $\texttt{p}(a_0,\dots,a_n)$ is a $k$-vector. The main reason to introduce free
$k$-terms is the following notion of equations.

\defi

A \emph{free} $k$\emph{-equation} (\emph{free equation} in short) is a pair
$\{\texttt{p}(x_0,\dots,x_n),\texttt{q}(x_0,\dots,x_{n'})\}$ of free $k$-terms. Given a fixed
equivalence relation $\sim$ on $\fin_k$, we will write the previous free equation as
\begin{equation*}
\texttt{p}(x_0,\dots,x_n)\rel \texttt{q}(x_0,\dots,x_{n'}).
\end{equation*}
Given $s,t$, $i_0$ and $i_1$-vectors respectively,  a free $j_0$-term $\texttt{p}$,  and
 a free $j_1$-term $\texttt{q}$ such that $\max\{i_l,j_l\}=k$ for $l=0,1$, we consider the equations of the
form $s+\texttt{p}\rel t+\texttt{q}$ and $\texttt{p}+s\rel \texttt{q}+t$, called
$k$\emph{-equations} (or \emph{equations}, if there is no possible confusion). The substitutions of
$(b_0,\dots,b_n)$ in the equation $s+p\rel t +q $ will be allowed only when $b_0>s,t$, and for an
equation $\texttt{p}+s\rel \texttt{q}+t$, provided that $b_n<s,t$. This last condition implies that
only finitely many substitutions are allowed for this latter equations, in contrast with the
equations of the form $s+p\sim t+q$.

\fdefi

\defi
We say that a $k$-equation $s+p(x_0,\dots,x_n)\rel t+q(x_0,\dots,x_n)$ (or
$p(x_0,\dots,x_n)\linebreak +s\rel q(x_0, \dots,x_n)+t$) \emph{holds} (or is \emph{true}) in $A$
iff for every $(a_0,\dots,a_n)$ in $[A]^{[n+1]}$ with $a_0> s,t$ (resp. $a_n<s,t$),
$s+p(a_0,\dots,a_n)\rel s+q(a_0,\dots,a_n)$ (resp. $p(a_0,\dots,a_n)+s\rel q(a_0,\dots,a_n)+s$).
The equation $s+p(x_0,\dots,x_n)\rel t+q(x_0,\dots,x_n)$ (or $p(x_0,\dots,x_n)+s\rel
q(x_0,\dots,x_n)+t$) is false in $A$ iff for every $(a_0,\dots,a_n)$ in $[A]^{[n+1]}$  with $a_0>
s,t$ (resp. $a_n<s,t$), $s+p(a_0,\dots,a_n)\nrel s+q(a_0,\dots,a_n)$ (resp.
$p(a_0,\dots,a_n)+s\nrel q(a_0,\dots,a_n)+s$). The equation is \emph{decided in $A$} iff it is
either true in $A$ or false in $A$.

It is clear that, given a $k$-equation $\texttt{p}(x_0,\dots,x_n)\rel
\texttt{q}(x_0,\dots,x_{n'})$, we can assume that $n=n'$, since we can extend the terms of the
equation adding summands of the form $T^k x$ and not changing the ``meaning" of the $k$-equation.

Some properties of equations that will be useful are given in the following.
\fdefi
\prop\label{preq}
Suppose that all free $k$-equations  with at most five variables  are decided in a given $k$-block
sequence $A$. Then:
\begin{enumerate}
\item If $x_0+T^{k-i} x_1 +x_2\rel x_0+ x_2$ is true in $A$, then   $x_0+T^{k-j} x_1 +x_2\rel x_0+
x_2$ is true in $A$ for every $j\le i$. \item If $x_0+x_1+Tx_2\rel x_0+Tx_2$  or $Tx_0+x_1+x_2\rel
Tx_0+x_2$ are true in $A$,  then  $x_0+x_1 +x_2\rel x_0+ x_2$ is also true in $A$. \item If the
equation $x_0+x_1+T^ix_2\rel x_0+T^ix_2$ is true in $A$, then the equation $x_0+x_1+T^jx_2\rel
x_0+T^jx_2$ also is true in $A$ for every $j\le i$.
 \item If the equation $T^{i}x_0+x_1+x_2\rel
T^ix_0+x_2$ is true in $A$, then  the equation $T^jx_0+x_1+x_2\rel T^jx_0+x_2$ also is true in $A$
for every $j\le i$. \item If the equation $x_0+T^{k-r_1}x_1+T^{k-r_0}x_2\rel x_0+T^{k-r_0}x_2$
holds, then also the equation $x_0+T^{k-r_2}x_1+T^{k-r_0}x_2\rel x_0+T^{k-r_0}x_2$ for every
$r_1>r_2$ and $r_0$.
\end{enumerate}

\fprop
\prue Suppose that the $k$-block sequence $A$ decides all the equations with at most five variables.

\noindent (i): Fix $j<i$. Then,
\begin{equation}
x_0+T^{k-i}x_1+T^{k-j}x_2+x_3 \sim x_0+T^{k-i}(x_1+T^{i-j}x_2)+x_3 \rel x_0+x_3 \text{ hold in
$A$.}
\end{equation}
Hence,
\begin{equation}
x_0+T^{k-i}x_1+(T^{k-j}x_2+x_3) \rel x_0+(T^{k-j}x_2+x_3) \text{ holds in $A$},
\end{equation}
and we are done.

\noindent (ii): Suppose now that $x_0+x_1+Tx_2\rel x_0+Tx_2$ is true in $A$. Then
\begin{equation}
\text{$x_0+x_2+Tx_3\rel x_0+Tx_3$ and $x_0+x_1+x_2+Tx_3\rel x_0+Tx_3$ are true in $A$.}
\end{equation}
Hence, $x_0+x_1+x_2+Tx_3\rel x_0+x_2+Tx_3$ holds in $A$, and therefore, $x_0+x_1+x_2\rel x_0+x_2$
is true in $A$.

\noindent (iii): Suppose that $x_0+x_1+T^ix_2\rel x_0+T^ix_2$ is true in $A$, and fix $j\ge i$.
Then, $x_0+x_1+x_2+T^j(x_3+T^{i-j}x_4)\rel x_0+x_1+x_2+T^jx_3+T^{i}x_4\rel x_0+T^{i}x_4$ hold in
$A$, and
\begin{equation}
\text{$x_0+x_1+T^j(x_2+T^{i-j}x_3)\rel x_0+x_1+T^jx_2+T^{i}x_3\rel x_0+T^{i}x_3$ hold in $A$,}
\end{equation}
which implies what we wanted.

\noindent (iv): This is showed in a similar manner that (iii).

\noindent (v):  Fix $r_1>r_2$ and $r_0$ and suppose that the equation
$x_0+T^{k-r_1}x_1+T^{k-r_0}x_2\rel x_0+T^{k-r_0}x_2$ holds in $A$. Then,
$x_0+T^{k-r_2}x_1+T^{k-r_1}x_2+T^{k-r_0}x_3\rel x_0+T^{k-r_1}(T^{r_1-r_2}x_1+x_2)+T^{k-r_0}x_3\rel
x_0+T^{k-r_0}x_3$ and $(x_0+T^{k-r_2}x_1)+T^{k-r_1}x_2+T^{k-r_0}x_3\rel x_0+T^{k-r_2}x_1\rel
T^{k-r_0}x_3$ holds in $A$. Therefore, $x_0+T^{k-r_2}x_1\rel T^{k-r_0}x_3\rel x_0+T^{k-r_0}x_3 $ is
true in $A$.
\fprue

\subsection{Systems of staircases,  canonical and staircase equivalence relations}
Classifying equivalence relations of $\fin_k$ is roughly the same as finding   properties of a
typical $k$-vector. One of these properties can be the cardinality, or, for example, the minimum or
maximum of its support.  Indeed Taylor's result on $\fin$  tells    that these are the relevant
properties of $1$-vectors. For an arbitrary $k>1$, one expects a longer list of   properties. One
example is obtained by considering for a given $k$-vector $a$ the least integer  $n$  of the
support of $a$ such that $a(n)=k$; another one is obtained by fixing  $i$ with $1\le i \le k$ and
considering the least $n$ such that $a(n)=i$. This is not always well defined, since for $i<k$
there are $k$-vectors where $i$ does not appear in their range. Nevertheless, this last property
seems very natural to consider. Indeed we  are going introduce a type of $k$-block sequences,
called systems of staircases, where these properties, and some others, are well defined for every
$k$-vector of their combinatorial subspaces.

\defi Given an integer $i\in [1,k]$ let $\mini_i,\maxi_i:FIN_k\to \N$ be the
mappings $\mini_i(s)=\min s^{-1}\{i\}$, $\maxi_i(s)=\max s^{-1}\{i\}$, if  defined, and 0
otherwise. A $k$-vector $a$ is a \emph{{system of staircases}} (\emph{sos} in short) \iff
\begin{enumerate}
\item $\rang s=\{0,1\dots,k\}$, \item $\min_ia<\min_j a<\max_j a<\max_i a$, for $i<j \le k$, \item
for every $1\le i\le k$,
\begin{align*}
 & \rang a \rest [\mini_{i-1}a,\mini_i a
]=\{0,\dots,i\}, \\
& \rang a \rest [\maxi_{i}a,\maxi_{i-1} a ]=\{0,\dots,i\},\\
& \rang a \rest [\mini_k a,\maxi_k a ]=\{0,\dots,k\}.
\end{align*}

\end{enumerate}
The following figure illustrates the previous definition.
\begin{center}
\begin{figure}[h]
\includegraphics[scale=0.8]{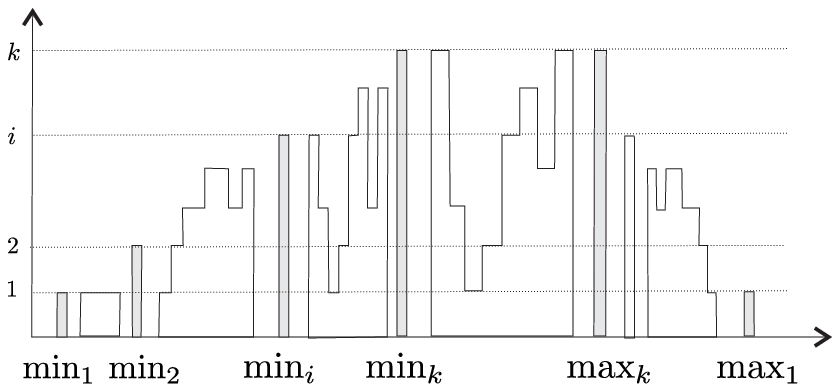}
\caption{ A typical sos.}
\end{figure}
\end{center}
A block subspace $A=(a_n)_n$  is a \emph{system of staircases} iff every $k$-vector in $  \langle A
\rangle$ is an sos. In the next proposition we show, among other properties, that for every
$k$-block sequence $A$ there is sos $B\in [A]^{[\infty]}$.
\fdefi

\prop\label{www}
\begin{enumerate}
\item[]
\item $T$ preserves sos, i.e., if $a$ is an sos $k$-vector, then $Ta$ is an sos $(k-1)$-vector.
\item $T^{k-j} a+b$, $a+T^{k-j}b$ are sos's, provided that $a<b$ are sos's. Therefore, for every $k$-term
$p(x_0,\dots,x_n)$ and every block sequence of sos $(a_0,\dots,a_n)\in [\fin_k]^{[n+1]}$, the
substitution $p(a_0,\dots,a_n)$ is also an sos.
\item A $k$-block sequence $A=(a_n)_n$ is an sos \iff  $a_n$ is an sos for
every $n$.
\item If $A$ is an sos, then any other $B\pe A$ is also  an sos.
 \item For every $A$ there is some $B\pe A$ which is an sos.
\end{enumerate}
\fprop
\prue

It is not difficult to prove   (i) and (ii) (for the last part of (ii), one can use induction on
the complexity of the $k$-term $p$). To show (iii), let us  suppose that    $a_n$ is an sos for
every $n$, and let us fix $a\in \langle (a_n)_n\rangle$. Then there is a $k$-term
$p(x_0,\dots,x_n)$ such that $p(a_0,\dots,a_n)=a$. Therefore, by (ii), $a$ is an sos. Assertion
(iv) easily follows from (ii). Finally, Let us prove (v): Fix $A=(a_n)_n$. For each $n$, let
\begin{equation*}
c_n=\sum_{j=1}^{k}T^{k-j}a_{(2k-1) n +j-1}+\sum_{j=1}^{k-1}T^{k-(k-j)}a_{(2k-1) n + k-1+j}.
\end{equation*}
Notice that for every $n$ one has that
$$\rang  c_n\rest[0,\text{min}_k (c_n)]=\rang
c_n\rest[\text{max}_k(c_n),\infty)=\{0,\dots,k\}.$$
 Therefore, $\rang T^{k-j} c_n\rest[0,\min_j
T^{k-j}(c_n)]=\rang T^{k-j}c_n\rest[\max_jT^{k-j}(c_n),\infty)=\{0,\dots,j\}$ for each $j\le k$.
 For $n\ge 0$, let
\begin{equation*}
b_{n}=\sum_{j=1}^{k}T^{k-j}c_{n(3k-1)+j-1}+
\sum_{j=1}^{k}T^{k-j}c_{n(3k-1)+k-1+j}+\sum_{j=1}^{k-1}T^{k-(k-j)}c_{n(3k-1)+2k-1+j}.
\end{equation*}
Now it is not difficult  to prove that every $b_n$ is an sos. \fprue
\defi
An equivalence relation $\rel$ on $FIN_k$ is \emph{canonical\footnote{this name is not arbitrary
chosen: We will show that every equivalence relation is, when restricted to some combinatorial
subpace, canonical.} in }$A$ \iff every $k$-equation are decided in every sos $B\in [A]$ in the
same way, i.e., iff for every $k$-equation $p\rel q$, either for every sos $B\in [A]$ one has that
 $p\rel q$ is true in $B$, or for every sos $B\in [A]$ one has that $p\rel q$ is false in $B$. We will say that
$\rel$ is {canonical} if it is canonical in $\fin_k$.
\fdefi

Canonical equivalence relations are those for which all the equations $p\sim q$ are decided in
every sos in the same way. It is not difficult to see that all the equivalence relations of the
list $\{\min,\max,(\min,\max), =, \fin^2\}$ are canonical in $\fin$.  Taylor's result for $\fin$
says that  there are no more canonical equivalence relations than the ones in this list. It will be
shown later that for every $k$ there is also a  finite list of canonical equivalence relations.
Indeed we will give an explicit description of how canonical equivalence relations look like.

In order to do the same to the equivalence relations in  $\fin_k$ we have to give a list of
relations naturally defined for a typical sos.
\defi\label{introdstair}
For a set $X$, a $k$-block sequence $A$, and an arbitrary map $f:\langle A \rangle\to X$ we define
the relation $R_f$ on $\langle A \rangle$ by $s R_f t$ \iff $f(s)=f(t)$. Whenever there is no
possible  confusion, we are going to use the notation $s f t$ instead of $sR_f t$. Now fix  an sos
$A$. Recall that $\mini_i(s)=\min\conj{n}{s(n)=i}$ for a given integer $i\in [1,k]$ and $s\in
\langle A\rangle$. This mapping can be interpreted as $\mini_i:\langle A \rangle \to FIN_{i}$ in
the following way
\begin{equation*}
\mini_i(s)(n)=\left\{\begin{array}{ll}i & \text{if }n=\mini_i(s)\\
0 &\text{otherwise}.
\end{array}
\right.
\end{equation*}
Extending this,  define, for $I\con \{1,\dots,k\}$,  the mapping $\mini_I:\langle A \rangle\to
\fin_{\max I}\con \fin_{\le k}$ by $\mini_I(s)(n)=i$   if $n=\mini_i(s)$, for $i\in I$ and 0
otherwise, i.e.,
  $\mini_I(s)=\conj{(\mini_i(s),i)}{i\in I}$, and extended by 0 in the rest. Similarly, let
\begin{equation*}
\maxi_i(s)(n)=\left\{\begin{array}{ll}i & \text{if }n=\maxi_i(s)\\
0 &\text{otherwise},
\end{array}
\right.
\end{equation*}
and let $\maxi_I:FIN_k\to FIN_{\max I}$ be defined by $\maxi_I(s)=\conj{(\maxi_i(s),i)}{i\in I}$,
again extended by 0. Clearly  $\mini_I=\bigvee_{i\in I}\mini_i$ and $\maxi_I=\bigvee_{i\in
I}\maxi_i$, where for two mappings $f,g:\langle A \rangle \to \fin_{\le k}$ we define $(f\vee
g)(s)=f(s)\vee g(s)$.

We  now introduce a more sophisticated   class of functions. For  $l\le i-1$, let
$\theta^0_{i,l},\theta^1_{i,l}:\langle A \rangle\to FIN_{l}$ be the mappings defined by
\begin{align*}
\theta_{i,l}^0(s)= & \conj{(n,l)} { n \in (\mini_{i-1}(s),\mini_{i}(s)) \, \& \, s(n)=l }, \text{ extended by
0, and}\\
\theta_{i,l}^1(s)=& \conj{ (n,l)}{n\in (\maxi_{i}(s),\maxi_{i-1}(s))\, \& \,s(n)=l }, \text{ extended by 0}.
\end{align*}
In other words,  for a given integer $n$
\begin{align*}
\theta^0_{i,l}(s)(n)= &\left\{\begin{array}{ll} l & \text{if $n\in (\mini_{i-1}(s),\mini_{i}(s))$ and $s(n)=l $} \\
0&\text{  otherwise, and}
\end{array}
\right. \\
\theta^1_{i,l}(s)(n)= &\left\{\begin{array}{ll} l & \text{if $n\in (\maxi_{i}(s),\maxi_{i-1}(s))$ and $s(n)=l$} \\
0&\text{  otherwise.}
\end{array}
\right.
\end{align*}
For example, for $k=4$, $i=3$, $l=2$ and a given  sos $4$-vector $s$, $\theta_{3,2}^2(s)$ is the $2$-vector
such that $(\theta_{3,2}^2(s))(n)=2$ for every $n$  such that

\noindent  \emph{(a)} $s(n)=2$, and

\noindent  \emph{(b)} $n$ is  in the interval between $\mini_{2}(s)$ (i.e., the first $m$ such that
$s(m)=2$) and $\mini_3(s)$ (i.e., the first $m$ such that $s(m)=3$), and it is zero otherwise.

 For
$1\le l\le k$, let
\begin{equation*}
\theta^2_{l}(s)=\conj{(n,l) } {n\in (\mini_k(s),\maxi_k(s))\, \& \,s(n)=l }\text{ extended by
zero.}
\end{equation*}
We illustrate this   with another example: For $k=4$, $l=3$ and  an sos $4$-vector $s$,
$\theta^2_3(s)$ is the $3$-vector with value $l=3$ in every element $n$ of the support of $s$ such
that

\noindent \emph{(a)} $s(n)=3$, and

\noindent  \emph{(b)}  $n$ is in between $\mini_4(s)$ and $\maxi_4(s)$, and 0 otherwise.

By technical convenience,  we declare $\theta_{i,-1}^{0}, \theta_{i,-1}^{1}$ and $\theta_{-1}^{2}$
 the $0$ mapping (hence, the equivalence relations associated to them  are all equal to $\fin_k^2$).
Also, for $I=\buit$, the mappings $\mini_{I}$ and $ \maxi_I$ are simply the $0$ functions, i.e.,
$0(s)(n)=0$ for all $s$ and $n$.

\fdefi
\nota\label{extens}
\noindent (i) Sometimes we will use $\min_i$ or $ \max_i$ as a integers instead of $i$-vectors,
i.e., for example $\min_i(s)$ will denote the unique integer $n$ such that $\min_i(s)(n)=i$.

\noindent (ii) Also, we can extend  the  mappings $f$ defined before for $\fin_k$ to all $\fin_{\le
k}$ by setting $\bar{f}(s)=f(s)$, if it is well defined, and $\bar{f}(s)=0$, if not. For example,
for a $(\le \! k)$-vector $s$, $\min_i(s)(n)=i$ iff  $i\in \rang s $ and $n$ is the minimum $m$
such that $s(m)=i$, and $\min_i(s)=0$ otherwise; and $\theta^0_{i,l}(s)$ will have the same
definition, provided that the mappings $\min_{i-1}$ and $\min_i$ are well defined for $s$, and so
on.

\fnota
\prop\label{staireasy}
Suppose that $l$ is such that $-1<l\le i-1$. Then,
\begin{enumerate}
\item $\rel_{\theta^0_{i,l}}\con \rel_{
 \mini_{i-1}}\cap\rel_{
\mini_i}$, $\rel_{\theta^1_{i,l}}\con \rel_{\maxi_{i}}\cap \rel_{\maxi_{i-1}}$, and $\rel_{\theta^2_l}\con
\rel_{\min_k} \cap \rel_{\max_k}$.
\item $\rel_{\theta^2_l}\con \rel_{\theta^2_{l+1}}$ and
 $\rel_{\theta^\vep_{i,l}}\con \rel_{\theta^\vep_{i,l+1}}$.
\end{enumerate}
\fprop \prue We  prove the result in (i) for $\theta^{0}_{i,l}$. The other cases can be shown in a similar
way. Suppose that $\theta^{0}_{i,l}(s)=\theta^{0}_{i,l}(t)$;  we  show that
$\min_{i-1}(s)=\min_{i-1}(t)$. Let $n$ be such that $\min_{i-1}(t)(n)=i-1$. By symmetry,   it
suffices to prove that $s(n)=i-1$. So, let $r$ be the unique integer such that $T^{k-C_A(t)(r)}a_r
(n)=i-1$. Note that $C_A(t)(r)\ge i-1$. There are two cases to consider:

\noindent \emph{(a)} $C_A(t)(r)= i-1$. Since $a_r$ is an sos, there is some $m\ge n$ such that
$T^{k-C_A(t)(r)}a_r (m)=l$, and hence $\theta^0_{i,l}(t)(m)=l$ and $\theta^0_{i,l}(s)(m)=l$. This
implies that $C_A(s)(r)=C_A(t)(r)$, and hence $T^{k-C_A(t)(r)}a_r\sqsubseteq s $. Hence,
$s(n)=T^{k-C_A(t)(r)}a_r (n)=i-1$.

\noindent \emph{(b)} $C_A(t)(r)>i-1$. Then, $\theta^0_{i,l}$ is well defined for
$T^{k-C_A(t)(r)}a_r$, and $\theta^0_{i,l}(T^{k-C_A(t)(r)}a_r)\sqsubseteq
\theta^0_{i,l}(t)=\theta^0_{i,l}(t)$, which implies that $T^{k-C_A(t)(r)}a_r\sqsubseteq s $, and
again we are done.

Let us now prove  the result for $\theta^2_l$ in  (ii). Suppose that $\theta^2_l(s)=
\theta^2_l(t)$, i.e.,
\begin{equation*}
\conj{n\in [\mini_k(s),\maxi_k(s)]}{s(n)=l}=\conj{n\in [\mini_k(s),\maxi_k(s)]}{t(n)=l}.
\end{equation*}
Let $n\in (\mini_k(s),\maxi_k(s))$ be such that $s(n)=l+1$. We  show that $t(n)=l+1$. Let $r$ be
the unique integer such that $T^{k-C_A(s)(r)} a_r(n)=l+1$. Then, $C_A(s)(r)\ge l+1$, and since
$a_r$ is an sos, $T^{k-C_A(s)(r)} a_r^{-1}\{l\}\neq \buit$. Moreover,

\clam $(T^{k-C_A(s)(r)} a_r)^{-1}\{l\}\cap (\mini_k(s),\maxi_k(s)) \neq \buit$.
\fclam
\prucl Let
$r_0,r_1$ be the unique integers such that $a_{r_0}(\mini_k(s))=a_{r_1}(\maxi_k(s))=k$. Observe that $r_0\le
r \le r_1$. There are two cases: If $r_0<r<r_1$, then we are done since $(T^{k-C_A(s)(r)} a_r)^{-1}\{l\}\cap
[\mini_k(s),\maxi_k(s)]=(T^{k-C_A(s)(r)} a_r)^{-1}\{l\}$ is non empty.

Suppose that $r_0=r$  (the case $r_1=r$ is similar).  Then, $C_A(s)(r)=k$, and $\min_k s=\min_k
a_r$. So, $(a_r)^{-1}\{l\}\cap (\mini_k(a_r),\maxi_k(s))\neq \buit$, since $a_r$ is an sos, and
therefore $\rang a_r \rest(\mini_k a_r,\max_k a_r)=\{0,\dots,k\}$. \fprucl Now that  for every
$m\in (T^{k-C_A(s)(r)} a_r )^{-1}\{l\} \cap (\min_k s, \max_k s) $ one has that $t(m)=l$, since
$(T^{k-C_A(s)(r)} a_r )^{-1}\{l\}\cap (\min_k s, \max_k s)\con \theta^2_l(t)$. By Proposition
\ref{basis}, $C_A(t)(r)=C_A(s)(r)$, and hence $T^{k-C_A(s)(r)} a_r \sqsubseteq t$, which implies
that $t(n)=T^{k-C_A(s)(r)} a_r (n)=s(n)=l$.

The second inclusion in (ii) is shown in a similar manner. The details are left to the reader.
\fprue

The collection  of mappings introduced in Definition \ref{introdstair}  can be divided into pieces
as follows.
\defi\label{defoff}
Let $\mathcal{F}_{\min}= \{\mathrm{\min}_1,\dots,\mathrm{\min}_k\}$, $
\mathcal{F}_{\max}=\{\mathrm{\max}_1,\dots,\mathrm{\max}_k\}$, $
\mathcal{F}_{\med^{\vep}}=\conj{\theta^\vep_{i,l}}{i\in \{1,\dots,k\}\,  \,l\in \{1,\dots,i-1\}
},\text{ for } \vep=0,1$, and  $ \mathcal{F}_\med=\conj{\theta^2_{l}}{l\in \{1,\dots,k\} }\cup
\{0\}$. Set $$\mathcal{F}=\mathcal{F}_{\min}\cup \mathcal{F}_{\max}\cup \mathcal{F}_{\med^0}\cup
\mathcal{F}_{\med^1}\cup \mathcal{F}_{\med}.$$ Given a $k$-block sequence $A$ we say that a
function $f:\langle A \rangle \to \fin_{\le k}$ is a \emph{staircase} function (in $A$) if it is in
the lattice closure of $\mathcal{F}$. An equivalence relation $\rel$ in $A$ is a \emph{staircase}
(in $A$) iff $\rel=\rel_f$ for some staircase mapping $f$.
\fdefi

\defi
Let  $f,g:\langle A \rangle\to\fin_k$ be two functions defined on the $k$-combinatorial subspace defined by   $A$.

\noindent (i) We say that $f$ and $g$ are \emph{incompatible}, and we write $f\perp g$, when $f(s)\perp f(s)$ for every
 $s \in \langle  A \rangle$.

\noindent (ii) We write $f<g$ to denote that $f(s)<g(s)$ for every $s\in \langle  A \rangle$.

\noindent (iii) We say that $f$ and $g $ are \emph{equivalent} (in $A$), and we write
$f\equiv g$, when $\rel_f\equiv \rel_g$, i.e., if $f$ and $g$ define the same equivalence relation
in $A$.
\fdefi
\nota
The family $\mc F$ is pairwise incompatible, i.e. if $f\neq g$ in $\mc F$ then $f \perp g$. Also, if $f<g$ then $f\perp g$.
\fnota

The following makes the notion of staircase relation more explicit.
\prop\label{stairare} Suppose that $A$ is an sos, and suppose that $f:\langle A \rangle\to \fin_{\le k}$. Then
 the following are equivalent:
\noindent (i) $f$ is staircase.

\noindent (ii) There are $I_\vep\con \{1,\dots,k\}$, $J_\vep\con\conj{j\in I_\vep}{j-1\in I_\vep}$,
$(l^{(\vep)}_j)_{j\in J_\vep}$  with $l^{(\vep)}_j\le j-1$ (for $\vep=0,1$)  and  $ l^{(2)}_k$
  such that
$$f= \mini_{I_0}\vee
\bigvee_{j\in J_0}{\theta^0_{j,l^{(0)}_j}}\vee {\theta^2_{l^{(2)}_k}} \vee
{\maxi_{I_1}} \vee \bigvee_{j\in J_1} {\theta^1_{j,l^{(1)}_j}}.$$ We say that
$(I_0,J_0,(l^{(0)}_j)_{j\in J_0},I_1,J_1,(l^{(1)}_j)_{j\in J_1},l^{(2)}_k)$ are the \emph{values}
of $f$.

\noindent (iii) Either $f=0$ or there is a unique sequence $f_0<f_1<\cdots <f_n$, $f_0\neq 0$  such that
$f\equiv\bigvee_{i=0}^n f_i$ in $A$.

\fprop
\prue
This decomposition  is a direct consequence of the fact that $\mc F$ is a pairwise incompatible
family and the inclusions exposed in Proposition \ref{staireasy}.
\fprue

\prop \label{stairr}  Fix  a staircase mapping $f$ with decomposition $f=f_0\cup\dots \cup f_n$ with $f_0<\dots<f_n$ in $\mc F$,
 an sos $A=(a_n)_n$ and $k$-vectors $s$ and $t$
of $A$. Then
\begin{enumerate}
\item $f(s)=f(t)$ \iff $f_i(s)=f_i(t)$ for every $0\le i \le
n$.
\item $f(s)=f(t)$ iff $f(s \rest\supp t)=f(t)$ and $f(t\rest \supp s)=f(s)$.\footnote{Notice
that $s\rest\supp t$ is not necessarily a $k$-vector, but we can still apply $f$ to it; see Remark
\ref{extens}.}
\item Suppose that $s_0,s_1<t_0,t_1$ are $(\le \! k)$-vectors of $A$ such that
$s_0+t_0$, $s_1+t_1$ and $s_0+t_1$ are $k$-vectors. If $f(s_0+t_0)=f(s_1+t_1)$, then
$f(s_0+t_0)=f(s_0+t_1)$.
\end{enumerate}
\fprop \prue   (ii) follows from the fact that $f_i<f_j$ for $i<j$. Let us check
(ii) using (i). We may assume  that $f\in \mathcal{F}$. There are several cases to consider.

\noindent  \emph{(a)}  $f=\min_i$. Suppose that $\min_i(s)=\min_i(t)$. Then, $i\in \rang s\rest
\supp t$, and hence $\min_i(s\rest \supp t)=\min_i s=\min_i t=\min_i(t\rest \supp s)$. Suppose now that
$\min_i s<\min_i t$. Then, $\min_i s<\min_i t \le \min_i (t\rest \supp s)$. So, $\min_i(t\rest \supp s)\neq
\min_i s$.

\noindent \emph{(b)}  $f=\max_i$ is shown in the same way.

\noindent \emph{(c)} $f=\theta^0_{i,l}$. Suppose that $\theta^0_{i,l}(s)=\theta^0_{i,l}(t)$. Then,
by \emph{(a)},
 $\min_{j}s=\min_j t\rest\supp s$, and  $\min_{j}t=\min_j s\rest\supp t$, where $j=i-1$ or $j=i$.
 Fix $n\in (\min_{i-1}(s),\min_i(s))$ such that $s(n)=l$.  Then, $t(n)=l$, and hence
 $\theta^0_{i,l}(t\rest s)(n)=l$. Now suppose that   $\theta^0_{i,l}(t\rest s)(n)=l$. Then,
 $t(n)=l$, and hence $s(n)=l$.

 Suppose that $\theta^0_{i,l}(s)=\theta^0_{i,l}(t\rest\supp s)$ and  $\theta^0_{i,l}(t)=\theta^0_{i,l}(s\rest\supp
 t)$. Then, $\min_j(s)=\min_j(t)$ for $j=i-1,i$. Fix $n$ such that $\theta^0_{i,l}(s)(n)=l$.
 Then, $\theta^0_{i,l}(t\rest\supp s)(n)=l$, which implies that $t(n)=l$.

\noindent \emph{(d)} The cases of $f=\theta^1_{i,l}$ and $f=\theta^2_l$ have a similar proof that
\emph{(c)}.

Let us prove (iii). To do this, fix $s_0,s_1,t_0,t_1$ as in the statement, and  suppose that
$f(s_0+t_0)=f(s_1+t_1)$. Suppose that $f=\min_i$. If $\min_i(s_0+t_0)=\min_i(s_0)$, then clearly
$\min_i(s_0+t_0)=\min_i(s_0+t_1)$. If not we have that  $\min_i(s_0+t_0)=\min_i(t_0)$, hence by our
assumptions $\min_i(s_1+t_1)=\min_i(t_0)$. Since $s_1<t_0$, it follows that
$\mini_i(s_1+t_1)=\mini_i(t_1)$ and we are done. Suppose now that $f=\max_i$.  If
$\max_i(s_0+t_0)=\max_i(s_0)$, then $\max_i(s_1+t_1)=\max_i(s_1)$ (now using the fact that
$t_1>s_0$), and therefore, $t_1$ is a $(<i)$-vector. So, $\max_i(s_0+t_0)=\max_i(s_0+t_1)$. If
$\max_i(s_0+t_0)=\max_i(t_0)$, then $\max_i(s_1+t_1)=\max_i(t_1)$ and we are done. Suppose now that
$f=\theta^0_{i,l}$ and suppose that $\theta^0_{i,l}(s_0+t_0)(n)=\theta^0_{i,l}(s_1+t_1)(n)=l$. If
$s_1(n)=l$, then $s_0(l)=l$,  and hence $(s_0+t_1)(n)=l$. If $t_1(n)=l$, then clearly
$(s_0+t_1)(n)=l$. By symmetry, we are done in this case.  The cases $f=\theta^0_{i,l}$ and
$f=\theta^2_l$  have a similar proof. We leave the details to the reader.
\fprue

\prop\label{st->can} Any staircase equivalence relation is canonical.
\fprop
\prue By Proposition
\ref{stairr}, it suffices to prove the result only for staircases functions  $f\in \mathcal{F}$.
So, we fix $f\in \mathcal{F}$, set $\rel=\rel_f$ and consider an equation $p(x_0,\dots,x_n)\rel
q(x_0,\dots,x_n)$ where $p(x_0,\dots,x_n)=\sum_{d=0}^{n}T^{k-m_d}x_d$ and
$q(x_0,\dots,x_n)=\sum_{d=0}^{n}T^{k-u_d}x_d$. Set $p^*=p\circ \texttt{x}^{-1}$ and $q^*=q\circ
\texttt{x}^{-1}$. So $p^*(d)=m_d$ and  $q^*(d)=u_d$ for $d\le n$ and 0 for the rest. Fix two sos's
$A$ and $B$ ($B$ can be equal to $A$), and suppose that $p(a_0,\dots,a_n)\rel_f q(a_0,\dots,a_n) $
for some $(a_0,\dots,a_n)\in [A]^{[n+1]}$. We  show that $p(b_0,\dots,b_n)\rel q(b_0,\dots,b_n)$
for every $(b_0,\dots,b_n)\in [B]^{[n+1]}$. There are several cases to consider depending on $f$.

\noindent \emph{(a) }$f=\min_i$.  Let $d_0$ be the first $d$ such that $m_d\ge i$, and $d_1$ be the
first $d$ such that $u_d\ge i$. Then $\min_i(p(a_0,\dots,a_n))=\min_i (T^{k-m_{d_0}})a_{d_0}$ and
$\min_i(q(a_0,\dots,a_n))=\min_i (T^{k-u_{d_1}})a_{d_1}$. Since $\min_i
(T^{k-m_{d_0}})a_{d_0}=\min_i (T^{k-u_{d_1}})a_{d_1}$, we have that $d_0=d_1$ (otherwise,
$a_{d_0}\perp a_{d_1}$). Hence $m_{d_0}=u_{d_1}$ (because $T^r a \perp T^s a$ if $r\neq s$). So $p$
and $q$ satisfy that for every $d<d_0$, both $m_d$ and $u_d$ are less than $i$ and
$m_{d_0}=u_{d_0}=i$. This implies that $\min_ip(b_0,\dots,b_n)=T^{k-m_{d_0}}b_{d_0}=\min_i
q(b_0,\dots,b_n)$.

\noindent\emph{(b)}  $f=\max_i$ has a similar proof.

\noindent \emph{(c)} $f=\theta^0_{i,l}$. By Proposition \ref{staireasy}, $\rel_{\theta^0_{i,l}}\con
\rel_{\min_{i-1}}\cap\rel_{\min_i}$. Hence $\min_{i-\vep}
p(a_0,\dots,a_n)=\min_{i-\vep}q(a_0,\dots,a_n)$ for $\vep=0,1$.
 Define, for $\vep=0,1$, $d_\vep$ as the least integer $d$  such that $p^*(d_j)=q^*(d_j)\ge i-1+\vep$. So,  $d_0\le
 d_1$ and
 \begin{align}
\theta^0_{i,l}p(a_0,\dots,a_n)=& \theta^0_{i,l}\sum_{d=d_0}^{d_1}
 T^{k-m_d}a_d \\
 \theta^0_{i,l}q(a_0,\dots,a_n)= & \theta^0_{i,l}\sum_{j=d_0}^{d_1}
 T^{k-u_d}a_d.
 \end{align}
 We see now that for every $d\in [d_0,d_1]$ either  $m_d$ and
 $u_d$ are both less than $l$ or $m_d=u_d$. To do this, suppose that $d\in [d_0,d_1]$ is such that $m_d\ge l$.
 Then $\theta^{0}_{i,l} T^{k-m_d}a_d
  \sqsubseteq \theta^{0}_{i,l}p(a_0,\dots,a_n)=\theta^{0}_{i,l}q(a_0,\dots,a_n)$. Since for $d\neq
  d''$ in $[d_0,d_1]$ one has that
    $T^{k-u_{d'}}a_{d'}\perp T^{k-m_{d}}a_{d}$, it follows that $T^{k-m_d}a_d\sqsubseteq
  T^{k-u_d}a_d$, and hence  $u_d=m_d$.

 \noindent  \emph{(d)} The cases $f=\theta^1_{i,l}$ and $f=\theta^2_{k}$ have a similar proof.
\fprue
Let us  now give some other properties of  equations for staircase equivalence relations.
\prop\label{red100} Suppose that $\rel$ is a staircase equivalence relation with values
$I_0,$ $J_0,$ $I_1,$ $J_1$, $(l^{(0)}_{j})_{j\in J_0}$,  $(l^{(1)}_{j})_{j\in J_1}$ and
$l^{(2)}_k$, and suppose that $A$ is an sos.
\begin{enumerate}
\item Let $0\le r_0<r_1\le r_2$. If $T^{k-r_0}x_0+T^{k-r_2}x_1+x_2\rel T^{k-r_0}x_0+x_2$ is true in $A$,
 then $r_1\notin I_0$.
\item  If $l^2_k=-1$,  then the equation $x_0+x_1+x_2\rel x_0+x_2$ is true in $A$. If $l^2_k\neq -1$, then
for every $0<l<l^2_k$ the equation $x_0+T^{k-l}x_1+x_2\rel x_0+x_2$ holds in $A$. \item Suppose
that  $i\notin I_0$, and let $j=\max I_0\cap [1,i]$.  Then the equation
$T^{k-j}x_0+T^{k-i}x_1+x_2\rel T^{k-j}x_0+x_2$ is true in $A$.
\item If $l^{(0)}_j=-1$, then the equation $T^{k-(j-1)}x_0+T^{k-(j-1)}x_1+x_2\rel
T^{k-(j-1)}x_0+x_2 $ is true in $A$.
\item Suppose that $l^{(0)}_j\neq -1$, and let $h<l^{(0)}_j$.
Then the equation $T^{k-(j-1)}x_0+ T^{k-h}x_1+x_2\rel T^{k-(j-1)}x_0+x_2$ is true in $A$.
\item Suppose that $p(x_0,\dots,x_n)$ is  a $(\le\! k)$-term, and suppose that $p(x_0,\dots,x_n)+
T^{k-l}x_{n+1} +x_{n+3}\rel p(x_0,\dots,x_n)+ T^{k-l}x_{n+2} +x_{n+3}$ holds in $A$. Then
$p(x_0,\dots,x_n)+ T^{k-l}x_{n+1} +x_{n+2}\rel p(x_0,\dots,x_n)+  x_{n+2}$ also holds.
\end{enumerate} The analogous symmetric results are also true.
\fprop
\prue  We give some of the proofs. The rest are quite similar, and the details are left to the reader.
The main idea is to use the decomposition  of $f=\bigvee_{i=0}^n f_i$ be the decomposition of $f$
into elements of $\mc F$ with $f_0<\dots<f_n$.

\noindent (i): Fix  $(a_0,a_1,a_2)\in [A]^{[3]}$. Then
$\min_{r_1}(T^{k-r_0}a_0+T^{k-r_2}a_1+a_2)=\min_{r_1}T^{k-r_2}a_1$, while
$\min_{r_1}(T^{k-r_0}a_0+a_2)=\min_{r_1}(a_2)$. Hence
$\min_{r_1}(T^{k-r_0}a_0+T^{k-r_2}a_1+a_2)\neq \min_{r_1}(T^{k-r_0}a_0+a_2)$.

For the rest of the points (ii) to (vi) one shows that in each case the corresponding equations for
$\sim_{f_i}$
 hold for every $0\le i\le r$, and then use Proposition \ref{stairr} to conclude  that the desired equation also holds.
\fprue
\defi
We call a staircase relation  a \emph{$\min$-relation} if its corresponding  set
$I_1=\emptyset$, and a \emph{$\max$-relation} if $I_0=\emptyset$.
\fdefi
\nota
\noindent (i) Proposition \ref{staireasy} states that if $l_{k}^{2}\neq -1$, then $k\in I_0\cap
I_1$. Hence if $\sim$ is a $\min$-relation or a $\max$-relation, then $l_k^2=-1$.

\noindent (ii) The equation $x+s\rel x+t$ is  true if $\rel$ is a $\min$-relation and the relation
$s+x\rel t+x$ is true if $\rel$ is a $\max$-relation.
\fnota

\section{The main Theorem}
The next theorem is the main result of this paper.

\teor\label{maint}
For every $k$ and every equivalence relation $\rel$ on $\fin_k$ there is an sos $B$ such that
$\rel$ restricted to $\langle B \rangle$ is a staircase equivalence relation.
\fteor

Again we  use  Taylor's result, now to expose the role of equations.
 Fix an equivalence relation $\sim$ on $\fin$. A diagonal procedure shows that we can find  a block
sequence $A=(a_n)_n$ such that for  every $i_0,i_1,i_2,i_3,j_0,j_1,j_2,j_3\in \{0,1\}$ and every
$s,t \in \langle A \rangle$, the equation
\begin{equation}
s+T^{i_0}x_0+T^{i_1}x_1+T^{i_2}x_2+T^{i_3}x_3\rel t+T^{j_0}x_0+T^{j_1}x_1+T^{j_2}x_2+T^{j_3}x_3
\text{ is decided in $A$.}
\end{equation}
For arbitrary $k$, the corresponding result is stated in Lemma  \ref{eqnde}.  We consider the same
cases considered in original Taylor's proof:

\noindent  \emph{(a)}\emph{ $x_0\rel x_1$ holds}. Then $\rel$ is $\langle A \rangle^2$ on $\langle
A \rangle$: Let $s,t\in \langle A \rangle$, pick $u>s,t$, and hence $s,t\rel u$.

\noindent \emph{(b)} \emph{$x_0\rel x_1$ is false, $x_0+x_1\rel x_0$ is true, and $x_0+x_1\rel x_1$
is false}. Let us check that $\rel$ is $\rel_{\min}$ on $\langle A \rangle$. Fix $s,t\in \langle A
\rangle$. Suppose that $s\rel_{\min}t$, and let $n$ be the least integer  such that $C_A(s)(n)=1$.
Then $s=a_n+s'$, $t=a_n+t'$, and using the fact that $x_0+x_1\rel x_0$ holds, $s,t\rel a_n$.
Suppose now that $s\nrel_{\min} t$, and suppose that $\min(s)<\min(t)$, and pick $n$ as before.
Then $s\rel a_n$, $a_n<t$, and $a_n\rel t$, a contradiction.

\noindent  \emph{(c)} \emph{$x_0\rel x_1$ is false, $x_0+x_1\rel x_0$ is false, and $x_0+x_1\rel
x_1$ is true}. Similar proof that 2. shows that $\rel$ is $\rel_{\min}$ on $\langle A \rangle$.

\noindent  \emph{(d)}  \emph{$x_0\rel x_1$ is false, $x_0+x_1\rel x_0$ and  $x_0+x_1\rel x_1$ are
false, and $x_0+x_1+x_2\rel x_0+x_2$ is true}. We show that $\rel$ is $\rel_{\min}\cap \rel_{\max}$
on $\langle A \rangle$. It is rather easy to prove that $\rel_{\min}\cap \rel_{\max}\con \rel$ on
$\langle A \rangle$. For the converse, suppose that $\max s\neq \max t$ and  $s\sim t$. We may
assume that  $\max s<\max t$. Let $n$ be the maximal integer $m$ such that $C_A(t)(m)=1$. Then,
$t=t'+a_n$, and hence the equation $s\rel t'+x_0\text{ holds}$ and hence $ \text{$t'+x_0+x_1\rel
t'+x_0$ also holds}$ which implies that $x_0+x_1\rel x_0$ holds, a contradiction. Notice that this
proves that if $x_0+x_1\rel x_0$ is false, then $\rel \con \rel_{\max}$. We assume that $\max
s=\max t$ but $\min s \neq \min t$. Suppose that $\min s<\min t$. We  show that $s\nrel t$. Suppose
again that $s\rel t$ and work for a contradiction. Let $n_0,n_1$ be the minimum and the maximum of
the support of $s$ in $A$ resp., and let $m_0$ be the minimum of the support of $t$ i $A$. Then
$s=a_{n_0}+s'+a_{n_1}$, $t=a_{m_0}+t'+a_{n_1}$. Using that the equation $x_0+x_1+x_2\rel x_0+x_2$
is true, we may assume that $s'=t'=0$. Since $n_0<m_1\le n_1$, either the equation $x_0+x_2\rel
x_1+x_2$ is true or the equation $x_0+x_1\rel x_1$ is true. But the first case implies that the
equations $x_0+x_3\rel x_1+x_2+x_3$ and $x_0+x_3\rel  x_2+x_3$ hold and hence $x_0\rel x_0+x_1$
holds, a contradiction.

\noindent  \emph{(e)}  \emph{$x_0\rel x_1$, $x_0+x_1\rel x_0$, $x_0+x_1\rel x_1$, $x_0+x_1+x_2\rel
x_0+x_2$ are false}. Then $\rel$ is $=$ on $\langle A \rangle$. Suppose that $s\rel t$, and suppose
that $s\neq t$. Since $x_0+x_1\rel x_0$ is false, then $\max s=\max t$ (see 4. above). Let $n$ be
the maximal integer $m<\max s$ such that $C_A(s)(m)\neq C_A(t)(m)$, and without loss of generality
 we assume  that $C_A(s)(n)=1$ and $C_A(s)(n)=0$. Then, $s=s'+a_n+s''$, and $t=t'+s''$, with
$t'<a_n$. Therefore the equation $s'+x_0+x_1\rel t''+ x_1$ holds, which implies that
$s'+x_0+x_1+x_2, s'+x_0+x_2\rel t''+x_2$ holds, and hence the equation $x_0+x_1+x_2\rel x_0+x_2$ is
true, a contradiction.

For arbitrary $k$ the proof is done by induction on $k$, making use of several lemmas. From now on
we fix an equivalence relation $\rel$ on $\fin_k$. Our approach is the following. By the pigeonhole
principle Theorem \ref{gow}, there is always an sos $A$ who decides a finite class of equations. It
turns out that two  kind of equations we are interested in are of the form $x_0+s\sim x_0+t$,
$s+x_0\rel t+x_0$ where $s$ and $t$ are $(k-1)$-vectors. The reason is that if they are decided,
then we can define naturally the $(k-1)$-equivalence relations
\begin{align*}
s\sim_0 t&\text{ iff } s+x_0\sim t+x_0\text{ holds,}\\
s\sim_1 t&\text{ iff } x_0+s\sim x_0+t\text{ holds.}
\end{align*}
and then use the inductive hypothesis to detect both $\sim_0$ and $\sim_1$  as $(k-1)$-staircase
equivalence relations. The next thing to do is to interpret   $\sim_0$ and $\sim_1$ as
$k$-relations $\sim_0'$ and $\sim_1'$, and then prove that in a suitable restriction $\sim\con
\sim_0'\cap \sim_1'$. Finally, a few more equations decided in some sos  will force the
decomposition $\sim=\sim_0'\cap \sim_1'\cap R$ for a suitable staircase relation $R$.

\lema\label{eqnde} There is some sos $A=(a_n)_n$ such that for every $5$-tuples
$\vec{i},\vec{j}\in \{0,\dots,k\}^5$, and every $(\le\!k)$-vectors $s$ and $t$ of $\langle
A\rangle$, the $k$-equation
$$s+\sum_{l=0}^4 T^{\vec{i}(l)}x_l\rel
t+\sum_{l=0}^4 T^{\vec{j}(l)}x_l\text{ is decided in $A$.}$$
\flema
\prue
We find a fusion sequence $(A_r)_r$ of $k$-block sequences, $A_r=(a^r_n)_n$ such that for every integer $r$
 the equations $s+\sum_{l=0}^4 T^{\vec{i}(l)}x_l\rel
t+\sum_{l=0}^4 T^{\vec{j}(l)}x_l$ are decided in $A_r$ for every $(\le\! k)$-vectors $s,t$ of
$\langle (a_i^i)_{i<r} \rangle$ and every $\vec{i},\vec{j}\in \{0,\dots,k\}^5$. Once we have done
this, the fusion sequence $A=(a_r^r)_r$ works for our purposes:
 Fix an equation $e$, $s+\sum_{l=0}^4 T^{\vec{i}(l)}x_l\rel
t+\sum_{l=0}^4 T^{\vec{j}(l)}x_l$, and let $r$ be the least integer such that $s,t$ are $(\le\!
k)$-vectors of $\langle(a_i^i)_{i<r}\rangle$. Then $e$ is decided in $A_r$, hence  it is also
decided in $A$.

We justify the existence of the demanded fusion sequence. Suppose we have already defined
$A_r=(a_n^{r})_n$.  Let $\mc L$ be the set of all the $k$-equations
 of the form
$$s+\sum_{l=0}^4 T^{\vec{i}(l)}x_l\rel t+\sum_{l=0}^4
T^{\vec{j}(l)}x_l$$
  where $s$ and $t$ are $(\le \!k)$-vectors in $\langle (a_i^i)_{i\le
r}\rangle_{\le k}$ and  $\vec{i},\vec{j}\in \{0,\dots,k\}^5$. Let
$$\La:[(a_n^{r})_{n\ge 1}]^{[5]}\to \{0,1\}^{\mc L}$$
 be
the finite coloring defined for each $(c_0,\dots,c_4)\in [(a_n^{r})_{n\ge 1}]^{[5]}$  and each
equation $e$ of the form $s+\sum_{l=0}^4 T^{\vec{i}(l)}x_l\rel t+\sum_{l=0}^4 T^{\vec{j}(l)}x_l \in
\mc L$ by $\La(c_0,\dots,c_4)(e)=0$ iff
$$s+\sum_{l=0}^4 T^{\vec{i}(l)}c_l\rel t+\sum_{l=0}^4
T^{\vec{j}(l)}c_l.$$ By Lemma \ref{gowgen}, there is $A_{r+1}\in [(a_n^{(r)})_{n\ge 1}]^{[\infty]}$
such that $\La$ is constant on $[A_{r+1}]^{[5]}$, which is equivalent to
 all the equations considered above being decided in $A_{r+1}$.
\fprue

\subsection{The inductive step. The relations $\sim_0'$ and $\sim_1'$}  Suppose that Theorem \ref{maint} holds for
$k-1$. Our intention is, of course, to prove the case for $k$.  To do this we first associate two
$k-1$-relations to our fixed $k$-relation $\rel$ as follows.
\lema \label{reg1} There is an
sos $A$ and  two staircase  $k-1$-equivalence relations
  $\rel_0$ and  $\rel_1$ on $\langle A \rangle_{k-1}$ such that for every $s,t\in \langle A \rangle_{k-1}$,
\begin{align}
\label{he4tuyorte}\text{the $k$-equation } s+x_0\rel t+x_0 & \text{ is true  in $A$ \iff  } s\rel_0 t,\text{ and} \\
\label{uiryy} \text{the $k$-equation }x_0+s\rel x_0+t &\text{ is true in $A$ \iff  } s\rel_1 t.
 \end{align}
Moreover $\rel_0$ and $\rel_1$ are such that for any two $(k-1)$-vectors  $s$ and $t$ of $A$,
\begin{align}
\label{uityjj} s\rel_0t & \text{ iff the $(k-1)$-equation $s+x\rel_0 t+x$ holds in $A$, and } \\
\label{uityjj1}s\rel_1t &\text{ iff the $(k-1)$-equation $x+s\rel_0 x+t$ holds in $A$.}
\end{align}
\flema
\prue
Let $B=(b_n)_n$ be an sos satisfying Lemma \ref{eqnde}. Then for $(k-1)$-vectors $s$ and $t$ of $B$
the $(k-1)$-equations $s+x_0\rel t+x_0$ are decided in $B$. Now define the  relation $\rel'$ on
$\langle B \rangle_{k-1}$ as follows. For $s,t\in \langle B \rangle_{k-1}$,
\begin{equation*}
\text{$s\rel' t$ iff $s+x_0\rel t+x_0$ holds in $B$.}
\end{equation*}
It is not difficult to see that $\rel'$ is an equivalence relation. By the inductive hypothesis
there is some $(k-1)$-block sequence $B'=(b_n')_n\in [(Tb_n)_n]^{[\infty]}$ and some canonical
equivalence relation $\rel_0$ such that $\rel'$ coincides with $\rel_0$ on $B'$ (since, by
Proposition \ref{st->can}, all staircase equivalence relations are canonical). The $k$-block
sequence $A=(Sb_n')_{n\ge 1}$ and the $k$-equivalence relation $\rel_0$ clearly satisfy
\eqref{he4tuyorte}. We prove assertion  \eqref{uityjj} for $\sim_0$.  To do this, suppose that
$s\rel_0 t$. Then the $k$-equation $s+x_0\rel t+x_0$ holds. Since the equation
$\text{$s+Tx_0+x_1\rel t+Tx_0+x_1 $ is decided}$, it must be true. It follows that for every
$k$-vector $b>s,t$  we have that $s+Tb\rel_0 t+Tb$. Since $\rel_0$ is canonical, we obtain that the
$(k-1)$-equation
\begin{equation}\label{jkuyyww}
\text{$s+x_0\rel_0 t+x_0$ holds in $A$,}
\end{equation}
as desired. Now assume that  (\ref{jkuyyww}) is true. Fix a $(k-1)$-vector $u>s,t$. Then
$s+u\rel_0t+u$, i.e., the $k$-equation $\text{$s+u+x_0\rel t+u+x_0$ holds}$. Hence
$\text{$s+x_0\rel t+x_0$ holds}$, that is $s\rel_0 t$.

We justify now the existence of a staircase $k-1$-equivalence relation $\rel_1$ and an sos $A$ such
that the statements  \eqref{uiryy} and \eqref{uityjj1} hold. We can find a fusion sequence
$(A_r)_r$, $A_r=(a_n^r)_n$, of $k$-block sequences of $A$, and a list $(\rel_a^n)_{a\in \langle
(a_i^i)_{i<r}\rangle_k}$ defined on $\langle A_r \rangle_{k-1}$ such that for every $s,t\in \langle
A_r\rangle_{k-1}$,
$$\text{$a+s\rel a+t$ \iff $s\rel_a^n t$.}$$
Let $A_\infty=(a_n^n)_n$ be the fusion sequence of $(A_r)_r$. Now for every $a\in \langle
A_\infty\rangle$ let $n(a)$ be unique integer unique $n$ such that $a\in \langle
(a_i^i)_{i<n}\rangle\setminus \langle (a_i^i)_{i<n-1}\rangle$. Define the finite coloring
$$c:\langle A_\infty\rangle\to \text{canonical equivalence relations on }FIN_{k-1}$$ by
$c(a)=\rel_a^{n(a)}$. By Lemma \ref{gowgen} there is some $A\in [A_\infty]^{[\infty]}$ in which $c$
is constant, with value $\rel_1$. We check that $A$ and $\rel_1$ satisfy what we want. Fix $a\in
\langle A\rangle$ and  two $k-1$-vectors   $s,t$ of $A$ with $a<s,t$; then $a\in \langle
\theta_{n(a)}\rangle$ and $s,t$ are $k-1$-block sequences of $A_{n(a)}$. So, $a+s\rel a+t$ \iff
$s\rel_a^{n(a)}t$ \iff $s \rel_1 t$. I.e. $x_0+s \rel x_0+t$ holds  iff $s\rel_1 t$. Notice that in
particular all equations $x_0+s\rel x_0+t$ are decided in $A$.

Let us prove now the assertion \eqref{uityjj1}. To do this, fix two $(k-1)$-vectors $s,t$ of $A$.
If $s\rel_1 t$, then $x_0+s\rel x_0+ t$. Given a $(k-1)$-vector $u<s,t$, choose  a $k$-vector $a<u$
in $\langle(Sb_n')_{n\ge 0}\rangle$. Then $a+u+s\rel a+u+t$, and this implies that $u+s\rel_1 u+t$;
in other words, the $(k-1)$-equation $x_0+s \rel_1 x_0+t$ holds. Suppose now that the
$(k-1)$-equation $x_0+s \rel_1 x_0+t$ holds. Pick $(k-1)$-vector $u<s,t$. Then the $k$-equation
$x_0+u+s\rel x_0+u+t$ is true, and hence also $x_0+s\rel x_0+t$ holds (since this equation is
decided).

Finally, we justify the existence of the fusion sequence $(A_r)_r$. Suppose we have already defined
$A_r=(a_n^r)_n$ fulfilling its corresponding requirements.  For every $a\in \langle (a_i^i)_{i <
r}\rangle$, put $\rel_a^{n+1}=\rel_a^n$. For every $a\in \langle (a_i^i)_{i\le r}\rangle\setminus
\langle (a_i^i)_{i<r}\rangle $, let $R_a$ be the relation on $\langle(a_n^r)_{n\ge 1}
\rangle_{k-1}$ defined by
$$\text{$s R_a t$ \iff  $a+s\rel a+t$.}$$
By the inductive hypothesis, we can find some   $B\pe (a_n^r)_{n\ge 1} $ such that for every $a\in
\langle (a_i^i)_{i\le r}\rangle\setminus \langle (a_i^i)_{i<r}\rangle $ the relation $R_a$ is
staircase when restricted to  $B$. Then $A_{r+1}=B$ satisfies the requirements.
\fprue

Roughly speaking,  the assertions (\ref{uityjj}) and (\ref{uityjj1}) tell that the $(k-1)$-relation
$\sim_0$ does not depend on the part of a $(k-1)$-vector before $\min_{k-1}$ and that $\sim_1$ does
not depend on the part of a $(k-1)$-vector after $\max_{k-1}$.  Indeed \eqref{uityjj} and
\eqref{uityjj1} determine the form of $\rel_0$ and $\rel_1$. To express this mathematically we
introduce the following useful notation.
\defi
For $l\le k$, let $\max_k^{l}:\fin_k\to \fin_k$ be defined by
$$({\max}_k^l s)(n)=\left\{\begin{array}
{ll} s(n) & \text{if }n\le \max_k(s) \text{ and }s(n)\ge l,\\
0 &\text{otherwise.}
\end{array}\right.$$
In other words $\max_k^l$ is the staircase function with values $I_0=\{l,\dots,k\}$,
$J_0=\{l+1,\dots,k\}$, for every $j\in J_0$, $l_j^{(0)}=l$, $l^{(2)}_k=l$ and $I_1=\{k\}$.
Symmetrically, we can define $\min_k^{l}$ by $\min_k^l(n)=s(n)$ iff $n\ge \min_k s$ and $s(n)\ge
l$, and 0 otherwise.
\fdefi

\prop\label{wqkiol} Suppose that  $\Rel$ is a staircase relation, and suppose that $A$ is an sos.
The following are equivalent:
\begin{enumerate}
\item For every   $k$-vectors $s,t$ of $A$,  one has that $s\Rel t$ iff
$x+s\Rel x+t$ holds in $A$.
\item Either $\Rel$ is a
$\max$-relation or there is some $\max$-relation $\Rel'$ and some $l\in \{1,\dots,k\}$ such that
$\Rel=\Rel'\cap \max_k^l$.
\end{enumerate}
The analogous result for $s+x\Rel t+x$ is also true.
 \fprop
 \prue Fix a staircase relation $\Rel$
with values $I_\vep,J_\vep,(l^{(\vep)}_j)_{j\in J_\vep}$ ($\vep=0,1$) and $l^{(2)}_k$ such that for
every $k$-vectors $s,t$ one has that   $s\Rel t$ iff $x+s\Rel x+t$ holds. Suppose that $I_0\neq
\emptyset$, since otherwise $\Rel$ is a $\max$-relation. Let $l=\min I_0$. We  show that
$I_0=\{l,l+1,\dots,k\}$, $J_0=\{l+1,\dots,k\}$, for every $j\in J_0$, $l_j^{(0)}=l$, $l^{(2)}_k=l$
and $k\in I_1$.  First we show that $l^{(2)}_k\neq -1$. If not, the equation $x_0+x_1+x_2\Rel
x_0+x_2 $ is true and hence the equation $x_1+x_2\Rel x_2$ is true, which implies that $l\notin
I_0$, a contradiction. If $l^{(2)}_k>l$, then the equation {$x_0+T^{k-l}x_1+x_2\Rel x_0+x_2$ is
true} and hence the equation {$T^{k-l}x_1+x_2\Rel x_2$ is true}, which implies again that $l\notin
I_0$. If $l^{(2)}_k<l$, then the equation {$T^{k-l^{(2)}_k}x_0+x_1\Rel x_1$ holds} and hence the
equation
\begin{equation}
\text{$x_0+T^{k-l^{(2)}_k}x_1+x_2\Rel x_0+ x_2$ holds,}
\end{equation}
which contradicts the definition of $l^{(2)}_k$.

We now show  that $I_0=\{l,\dots,k\}$. It is clear that $I_0\con \{l,\dots,k\}$ since $l$ is the
minimum of $I_0$. We  prove the reverse inclusion $\{l,\dots,k\}\con I_0$. Suppose not, and set
\begin{equation*}
j=\min \{l,\dots,k\}\setminus I_0.
\end{equation*}
Then the equation {$T^{k-j-1}x_0+T^{k-j}x_1+x_2\Rel T^{k-j-1}x_0+x_2$ is true} and  hence the
equation {$x_0+T^{k-j-1}x_1+T^{k-j}x_2+x_3\Rel x_0+ T^{k-j-1}x_1+x_3$ is true}, which implies that
the equation $x_0+T^{k-j}x_1+x_2\Rel x_0+ x_2$ also holds. This  contradicts the fact that $j>l$
and that $\rel\con \Rel_{\theta^2_l}$. Notice that $I_0=\{l,\dots,k\}$ implies that
$J_1=\{l+1,\dots,k\}$.

We show that $l^{(0)}_j=l$  for all $j\ge l+1$. Suppose that $l^{(0)}_j=-1$. This implies that the
equation {$T^{k-(j-1)}x_0+T^{k-(j-1)}x_1+x_2\Rel T^{k-(j-1)}x_0+x_2$ holds}. Again by  adding  one
variable at the beginning of both terms and using the fact that  $j-1\ge l$ we can  arrive
 at a contradiction to  the fact that $l^{(2)}_k=l$.  Suppose now that $l_j^{(0)}<l$. Then the equation
$T^{k-l^{(0)}_j}x_0+x_1\Rel x_1$ is true, and adding a variable we arrive at a contradiction.
Suppose that $l_j^{(0)}>l$, then the equation
\begin{equation}
\text{$T^{k-(j-1)}x_0+T^{k-l}x_1+x_2\Rel T^{k-(j-1)}x_0+x_2$ is true,}
\end{equation}
which yields  a contradiction in the same way as before. It is not difficult to check that the
converse and the analogous  situation for $\min$ are also true.
\fprue
Proposition \ref{wqkiol} and (\ref{uityjj}) and (\ref{uityjj1})  determine the relations $\rel_0$
and $\rel_1$ as follows.
\cor The relation $\rel_0$ is either a $\min$-relation or there is some
$l\le k-1$ and some $\min$-relation $R$ such that $\rel_0=R\cap \min_{k-1}^l$ and $\rel_1$ is
either a $\max$-relation or there is some $l\le k-1$ and some $\max$-relation $R$ such that
$\rel_1=R\cap \max_{k-1}^{l}$. \qed\fcor

Recall that $\rel_0$ and $\rel_1$ are both staircase equivalence relations of $\fin_{k-1}$. We now
give  the proper interpretation of both  as $k$-relations. Suppose that $k>1$. We know that either
$\rel_1$ is a max-relation, or $\rel_1=\max_{k-1}^l\cap R$, with $R$ a $\max$-relation. Let

$$\rel_1' =\left\{\begin{array}{ll}
\rel_1 & \text{if } \rel_1 \text{ is a  $\max$-relation} \\
\theta^1_{k,l}\cap R & \text{if $I_0\neq \buit$.}
\end{array}\right.$$
Notice that in the second case we have that $\max_k\con \sim_1'$. We do the same for $\rel_0$: It
is either a $\min$-relation or $\rel_0=R\cap \min_{k-1}^l$, being $R$ a $\min$-relation. Let
$$\rel_0'=\left\{\begin{array}{ll}
\rel_0 & \text{if } \rel_0 \text{ is a  $\min$-relation} \\
R\cap \theta^0_{k,l} & \text{if $I_1\neq \buit$.}
\end{array}\right.$$
In this second case we have that $\min_k\con \sim_0'$. For $k=1$, let $\rel_0'=\rel_1'=\fin_1^2$.

So, although $\rel_0$ is not a $\min$-relation and $\rel_1$ is not a $\max$-relation, their
corresponding interpretations $\rel_0'$ and $\rel_1'$ as $k$-relations are a $\min$-relation and a
$\max$-relation, respectively.

The relations $\rel_0'$ and $\rel_1'$ have similar properties than  $\rel_0$ and $\rel_1$.

\prop\label{stairr2} Let $s$ and $t$ be $(k-1)$-vectors. Then
\begin{enumerate}
\item $s+x\rel_0' t+x$ holds iff $s\rel_0 t$ iff $s+x\rel t+x $ holds.
\item $x+s\rel_1' x+t$ holds iff $s\rel_1 t$ iff $x+s\rel x+t$ holds.
\item $s+x_0+x_1\sim_0' t+x_0+x_2$ holds iff $s \rel_0 t$.
\item $x_0+x_2+s\sim_1' x_1+x_2+t$ holds iff $s \rel_1 t$.
\end{enumerate}
\fprop
\prue
We show the result for $\rel_1$; for $\sim_0$ the proof  is similar, and we leave the details to
the reader. If $\rel_1$ is a $\max$ relation, then there is nothing to prove. Suppose that $\rel_1=
\max_{k-1}^l\cap R$, for some $l\le k-1$, where  $R$ is a $\max$ relation. So, $\rel_1'=
\theta^{(1)}_{k,l}\cap R$ and we only have to show that
\begin{equation}
\text{$s\,\mathrm{max}_{k-1}^l t$ iff the equation $x+s \, \theta^{(1)}_{k,l} x+t $ holds,}
\end{equation} which is not difficult  to check (see Figure below).
\begin{center}
\begin{figure}[h]
\includegraphics[scale=0.6]{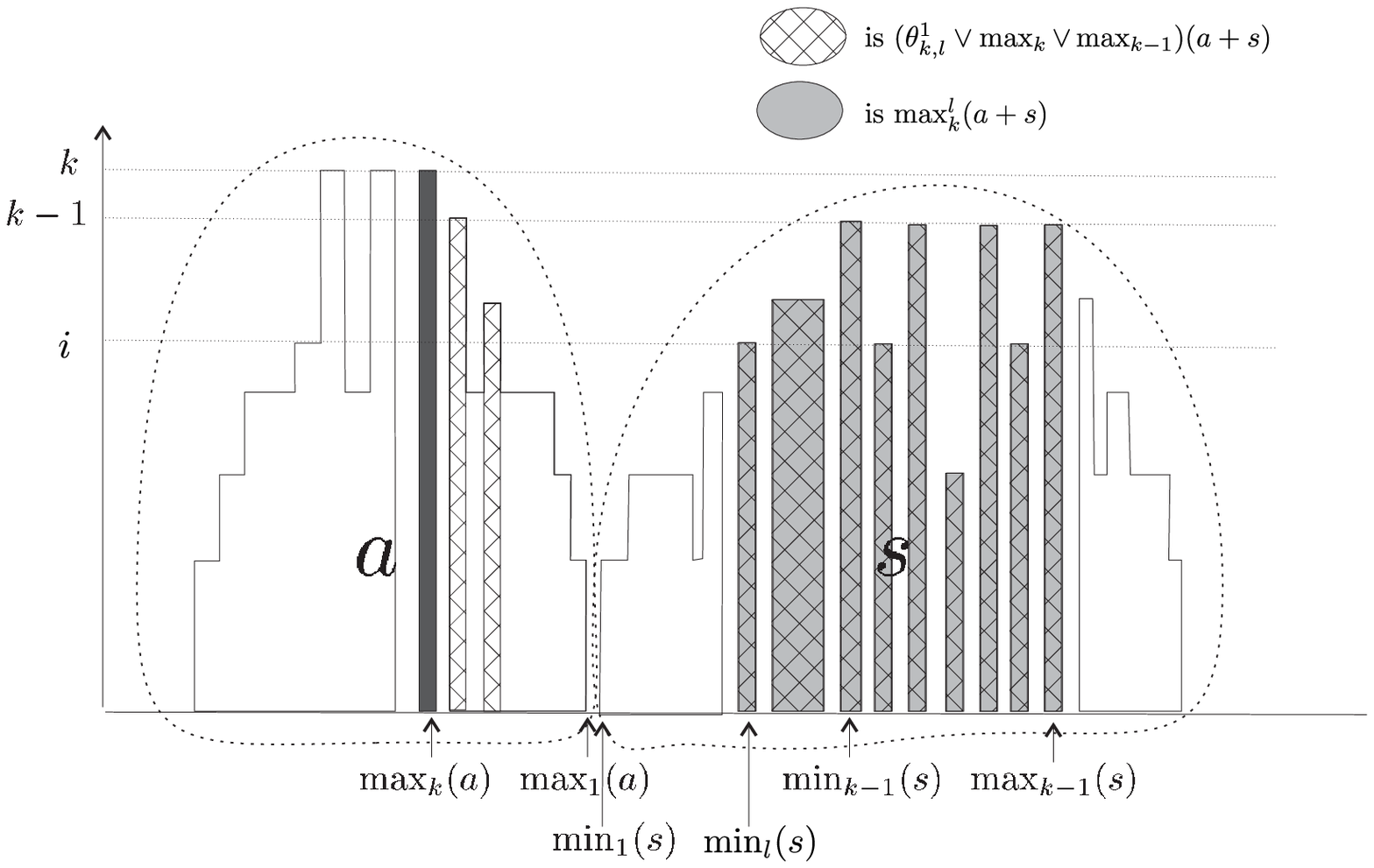}
\caption{ The relation between $\rel_1$ and $\rel_1'$}
\end{figure}
\end{center}
\fprue

\defi
Let $D=(d_n)_n$  be a $k$-block sequence and a $k$-vector $s=\sum_{n\ge 0}T^{k-C_D(s)(n)}d_n$ of
$D$, and let $n_0=n_0(s)$ and $ n_1=n_1(s)$ be respectively the minimal and the maximal elements of
the set of integers $n$ such that $C_D(s)(n)=k$. We define the \emph{first part of $s$ in $D$} as
the $(\le \! k-1)$-vector $f_D s=\sum_{n<n_0}T^{k-C_D(s)(n)}d_n$,  the \emph{middle part of $s$ in
$D$} as the $(\le \! k)$-vector $m_D s=\sum_{n\in (n_0,n_1)}T^{k-C_D(s)(n)}d_n$ and the \emph{last
part of $s$ in $D$}, as the $(\le \! k-1)$-vector $l_Ds=\sum_{n>n_1}T^{k-C_D(s)(n)}d_n$. Using
this, we have  the   decomposition
$$s=f_Ds+b_{n_0}+m_Ds+b_{n_1}+l_Ds.$$
So $f_D s$ is the part of $s$ before the occurrence of $\min_{k}s$, $m_D s$ is the part of $s$
between $\min_k s$ and $\max_k s$, and $l_D s$ is the  part of $s$ after $\max_{k}s$. All these
definitions are local, depending on a fixed sos $D$.

Let $\mathbb{A}=(a_n)_n$ satisfy both Lemmas \ref{eqnde} and \ref{reg1}, and let
$\mathbb{B}=(b_n)_n$ be defined for every $n$ by $b_n=T a_{3n}+a_{3n+1}+Ta_{3n+2}$.  The  role of
$\mathbb{B}$ is to guarantee that for every $k$-vector $s$ of $\mathbb{B}$ the first part
$f_\mathbb{A} s$ and the last part  $l_\mathbb{A}s$ are  both $(k-1)$-vectors. We need this because
$\rel_\vep$ ($\vep=0,1$) gives information only about $(k-1)$-vectors, since it is a
$k-1$-relation.

\fdefi

From now on we work in $\mathbb{B}$, unless we  explicitly say the contrary. The following
proposition tells us that many equations are decided in $\mathbb{B}$.
\prop\label{red88}
Let  $p(x_1,\dots,x_{n-1})$ and $q(x_1,\dots,x_{n-1})$ be  $(\le\!k-1)$-terms. Then:
\begin{enumerate}
\item The equation $x_0+p(x_1,\dots,x_{n-1})\rel x_0+q(x_1,\dots,x_{n-1})$ is decided in $\mathbb{B}$.
\item The equation $x_0+p(x_1,\dots,x_{n-1})\rel x_0+q(x_1,\dots,x_{n-1})$ holds in $\mathbb{B}$ iff
the equation $x_0+p(x_1,\dots,x_{n-1})\rel_1' x_0+q(x_1,\dots,x_{n-1})$ holds in $\mathbb{B}$.
\end{enumerate}
The analogous   results for $\rel_0'$ are also true.
\fprop
\prue
Fix two  $(\le\!k-1)$-terms $p=p(x_1,\dots,x_{n-1})$, $q=q(x_1,\dots,x_{n-1})$.

\noindent (i) Fix a finite block sequence $(c_0,\dots,c_{n-1})$ in $\mathbb{B}$. Suppose that
$c_0+p(c_1,\dots,c_{n-1})\rel c_0+q(c_1,\dots,c_{n-1})$. By definition of $\mathbb{B}$,
$c_0=c_0'+c_0''$, where $c_0'$ is a $k$-vector of $\mathbb{A}$ and $c_0''$ is a $(k-1)$-vector of
$\mathbb{A}$. Hence,
\begin{equation}
c_0''+p(c_1,\dots,c_{n-1})\rel_1 c_0'' +q(c_1,\dots,c_{n-1}).
\end{equation}
Since the relation $\rel_1$ is $(k-1)$-canonical in $\mathbb{A}$, the $(k-1)$-equation
\begin{equation}
\text{$x_0+p(x_1,\dots,x_{n-1})\rel_1 x_0+q(x_1,\dots,x_{n-1})$ is true in $\mathbb{A}$.}
\end{equation}
Fix $(d_0,\dots,d_{n-1})$ in $\mathbb{B}$, and set $d_0=d_0'+d_0''$. Then,
\begin{equation}
d_0''+p(d_1,\dots,d_{n-1})\rel_1 d_0''+q(d_1,\dots,d_{n-1}),
\end{equation}
and hence, the equation
\begin{equation}
\text{$x_0+p(d_1,\dots,d_{n-1})\rel x_0+q(d_1,\dots,d_{n-1})$ holds in $\mathbb{A}$,}
\end{equation}
which implies that $d_0+p(d_1,\dots,d_{n-1})\rel_1 d_0+q(d_1,\dots,d_{n-1})$, as desired.

\noindent (ii) Suppose that $x_0+p(x_1,\dots,x_{n-1})\rel x_0+q(x_0,\dots,x_{n-1})$ holds in
$\mathbb{B}$. Then for a given block sequence $(c_0,c_1,\dots,c_{n-1})$ in $\mathbb{B}$, the
equation
\begin{equation}\label{iorepe}
\text{$x_0+Tc_0+p(c_1,\dots,c_{n-1})\rel x_0+Tc_0+q(c_1,\dots,c_{n-1})$ holds in $\mathbb{B}$.}
\end{equation}
By Proposition \ref{stairr2},  the assertion (\ref{iorepe}) implies that
\begin{equation}
\text{$x_0+Tc_0+p(c_1,\dots,c_{n-1})\rel_1' x_0+Tc_0+q(c_1,\dots,c_{n-1})$ holds in $\mathbb{B}$.}
\end{equation}
Since $\rel_1'$ is canonical, the equation
\begin{equation}
\text{$x_0+Tx_1+p(x_2,\dots,x_{n})\rel_1' x_0+Tx_1+q(x_2,\dots,x_{n})$ holds in $\mathbb{B}$.}
\end{equation}
Therefore the equation $x_0+p(x_1,\dots,x_{n-1})\rel_1' x_0+q(x_1,\dots,x_{n-1})$ holds in
$\mathbb{B}$, as desired.
\fprue
\prop\label{red3}
 Suppose that $a,b$ are $k$-vectors of $\mathbb{B}$, $s,t$ are $\le\!(
k-1)$-vectors of $\mathbb{B}$ such that $a<s$, $b<t$ and suppose that $a+s \rel_1' b+t$.
\begin{enumerate}
\item If $a,b < s,t$, then $a+s\rel
a+t$.
\item  If $l_\mathbb{A}a=l_\mathbb{A}b=0$, and $\max_k(a)>\max_k(b)$, then $b+s\rel b+t $.
\end{enumerate}
The corresponding analogous  results for $\rel_0'$ are also true.
\fprop
\prue Let us check (i): By
point (iv) of Proposition \ref{stairr}, we have that $a+s \rel_1' a+t$. By construction of
$\mathbb{B}$, $a=a'+a''$ where $a'$ is a $k$-vector and $a''$ is a $(k-1)$-vector, both of
$\mathbb{A}$. But since the relation is $\rel_1'$ is staircase, it is canonical, and hence the
$k$-equation
\begin{equation}
\text{$x_0+a''+s \rel_1' x_0+a''+t$ holds in $\mathbb{A}$.}
\end{equation}It follows from Proposition \ref{red88}   that $a''+s \rel_1
a''+t$, and hence, by definition of $\sim_1$, the $k$-equation
\begin{equation}\label{nbfgla}
\text{$x_0+a''+s\rel x_0+a''+t$ holds in $\mathbb{A}$.}
\end{equation}
Replacing in (\ref{nbfgla}) $x_0$ by $a'$, we obtain that $a+s\rel a+t$.

\noindent (ii):  Since $l_\mathbb{A}a=l_\mathbb{A}b=0$, we have that $a+s=a'+a_{n_0}+s$ and
$b+t=b'+a_{m_0}+t$. Since $\max_k(a)>\max_k(b)$, it follows that $n_0>m_0$.  This together with the
fact that $a+s\rel_1' b+t$ implies that  $\max_k\not\con \rel_1'$ and hence, by definition,
$\sim_1'$ has to be  $\max$-relation. Set $i=\max I_1(\rel_1')<k$. Then the equation
\begin{equation}\label{eq789}
\text{$p(x_0,\dots,x_r)+T^{k-i'}x_{r+1}\rel_1' q(x_0,\dots,x_r)+T^{k-i'}x_{r+1}$ is true,}
\end{equation}  for every
terms $p$ and $q$, and every $i'\ge i$. Now set
\begin{equation*}
t=t'+T^{k-j}a_{n_0}+t''.
\end{equation*}
Notice that $t''$  is an $i$-vector, and $s$ is an $i'$-vector for some $i'\ge i$. By
(\ref{eq789}),
\begin{equation}
a+s=a'+a_{n_0}+s\rel_1' b'+a_{m_0}+t'+T^{k-j}a_{n_0}+s\rel_1' b'+a_{m_0}+s=b+s.
\end{equation}
Hence, $b+t\rel_1' b+s$, and since $b<s,t$,  1. implies that  $b+s \rel b+t$.
\fprue

Our intention is to show that $\rel\con \rel_1'$. To do this, we decompose the relation $\rel_1'$
as the final step of a chain $\rel_1'(1)\con \cdots \con \rel_1'(k)=\rel_1'$ and we prove by
induction on $j$ that $\rel\con\rel_1'(j)$.
\defi
Suppose that $\Rel$ is a $\max$-relation with values $I_0=\buit,I_1,J_1$ and $(l_j^{(1)})_{j\in
J_1}$. For every $i\le k-1$ we define  $I_1(i)=I_1\cap [0,i]$, $J_1(i)= J_1\cap [0,i]$, and let
 $\Rel(i)$ be the staircase equivalence relation on $\fin_k$ with values $I_1(i)$, $J_1(i)$, and $
(l_j^{(1)})_{j\in J_1(i)}$. So the relations between $\Rel(i+1)$ and $\Rel(i)$ is the following:
\begin{align*}
\Rel(i+1)=&\left\{ \begin{array}{ll} \Rel(i) & \text{if } i+1\notin I_1\\
\maxi_{i+1}\cap \Rel(i) & \text{if } i+1\in I_1 \text{ and } i\notin I_1\\
\maxi_{i+1}\cap \Rel(i)\cap \theta^1_{i+1,l^{(1)}_{i+1}} & \text{if } i+1\in J_1,
\end{array}\right.
\end{align*}
and  $\Rel=\Rel(k)$. Observe that each  $\Rel(i)$ is also a staircase equivalence relations on
every sos of $\fin_i$.
\fdefi
Roughly speaking, $\Rel(i)$ is the staircase equivalence relation whose values are the ones from
$\Rel$ which are smaller than $i$.
\nota
One has that for a given $i\le k-1$,  $s\Rel(i)t$ iff the equation with variable $x$
\begin{equation}
\text{$x+s\rest [\maxi_i(s),\maxi_1(s)]\Rel(i)x+ t\rest [\maxi_i(s),\maxi_1(s)]$ holds.}
\end{equation}
\fnota
\prop Suppose that $\Rel$ is a $\max$-relation of $\fin_k$.
Fix  $j'<j<j''$, and suppose that $s$ is a $j'$-vector, $t$ is a $(<j)$-vector, and $a$ is a
$j''$-vector such that $a+s\Rel(j) T^l a+ t$ for some $l>0$. Then, $\Rel(j)=\Rel(j')$, and hence
$s''\Rel(j) t''$.
\fprop
\prue
Set $s'=a+s$ and $t'=T^la+t$, and suppose that $s'\Rel(j)t'$. We are going to show that
$I_2(j)=I_2(j')$, which will imply that $\Rel(j)=\Rel(j')$, as desired. We know that
$s'\rest[\maxi_j(s'),\maxi_1(s')]\Rel(j)t'\rest[\maxi_j(s'),\maxi_1(s')]$. Notice that for every
$r\in [j,{j'})$, $\maxi_r(s')=\maxi_r(a)$, hence   $\maxi_r(s')\neq \maxi_r(t')$ , since $a$ and
$T^l a$ have nothing in common except 0's.   This implies that $I_2(j)\con [j',1]$ and hence
$I_2(j)=I_2(j')$.
\fprue

\lema\label{rel1}
$\rel\con\rel_1'(j)$, for every $j\le k$. In particular, $\rel\con \rel_1'$.
\flema
\prue
The proof is by induction on $j$. Notice that if $k=1$, then $\rel_1'=\fin_1^2$ and hence there is
nothing to prove. Suppose that $k>1$. Let $I_1$, $J_1$ and $(l^{(1)}_j)_{j\in J_1}$ be the values
of $\rel_1'$.

\underline{$j=1$}: Suppose that $1\in I_1$ (otherwise there is nothing to prove), i.e.,
$\rel_1'(1)=\rel_{\max_1}$. Suppose that $s\rel t$ but $\max_1(s)<\max_1(t)$, and let $n$ and $i$
be the unique integers such that
\begin{equation}
\text{$\maxi_1 T^{k-i}a_n=\maxi_1 t$ and $t=t'+T^{k-i}a_n$.}
\end{equation}  So, $s=s'+T^{k-i'} a_n$, for some $i'<i$ and some $k$-vector $s'$. The fact that $s\sim t$ implies that
the equation $s'+T^{k-i'}x_0 \rel t'+T^{k-i}x_0$ holds in $\mathbb{B}$, which implies that the
equation $s'+T^{k-i'}(x_0+T^{i'}x_1)\rel t'+T^{k-i}(x_0+T^{i'}x_1)$ is true. Therefore
\begin{equation}\label{hlkop}
\text{$s'+T^{k-i'}x_0+\rel t'+T^{k-i}x_0+T^{k-i+i'}x_1$ holds,}
\end{equation}
(\ref{hlkop}) implies that the equation
\begin{equation}
\text{$t'+T^{k-i}x_0\rel t'+T^{k-i}x_0 + T^{k-i+i'}x_1$ is true,}
\end{equation}
and hence, also
\begin{equation}
\text{$x_0+T^{k-i+i'}x_1\rel x_0$ is true in $\mathbb{B}$.}
\end{equation}
But since $j-i+i'<k$, we have that
\begin{equation}
\text{$x_0+Tx_1\rel x_0+Tx_1+T^{k-i+i'}x_2$ is true,}
\end{equation}
and by  Proposition \ref{red88}, we have that
\begin{equation}
\text{$x_0+Tx_1\rel_1' x_0+Tx_1+T^{k-i+i'}x_2$ holds,}
\end{equation}
which contradicts the fact that $1\in I_1$.

\noindent \underline{$j\curvearrowright j+1$}. Assume that  $\rel\con\rel_1'(j)$ and let us
conclude that $\rel\con\rel_1'(j+1)$. There are two cases:

\noindent {\emph{(a)} {$j\notin I_1$}}: Suppose that $j+1\in I_1$ (otherwise, there is nothing to
prove), and set
\begin{equation}
\be=\max I_1\cap [0,j].
\end{equation}
Notice that $\be$ can be 0. By definition of  $\rel_1'$, we know that  if $j+1=k$  belongs to
$I_1$, then $j=k-1$ also belongs to $I_1$. So, $j+1<k$. We only need to show that $\rel\con
\max_{j+1}$: Suppose that $s\rel t$, and $\max_{j+1}s<\max_{j+1}t$; set $s=s'+T^{k-l}a_n+s''$,
$t=t'+T^{k-l'}a_n +t''$, with $l<l'$, $l'\ge j+1$, and $(<(j+1))$-vectors  $s''$ and $t''$. Observe
that in the previous decomposition of $s$, $s'$ needs to be a $k$-vector. By the inductive
hypothesis,
\begin{equation}
s'+T^{k-l}a_n+s'' \rel_1'(j) t'+T^{k-l'}a_n +t''.
\end{equation} Since $\rel_1'$ is a staircase equivalence relation,
(iv) of Proposition \ref{stairr} gives  that
\begin{equation}\label{wlky}
\text{$s'+T^{k-l}a_n+s'' \rel_1'(j) s'+T^{k-l}a_n +t''$ ($t'' $ can be 0)},
\end{equation}
which implies that $s'+T^{k-l}a_n+s'' \rel_1' s'+T^{k-l}a_n +t''$, and hence, by Proposition
\ref{red3}, $s'+T^{k-l}a_n+s'' \rel s'+T^{k-l}a_n +t''$. Resuming, we have that
\begin{equation}
s'+T^{k-l}a_n+t''\rel t'+T^{k-l'}a_n+t'',
\end{equation}
and hence, the equation
\begin{equation}
\text{$s'+T^{k-l}x_0+T^{k-\al}x_1\rel t'+T^{k-l'}x_0+T^{k-\al}x_1$ holds,}
\end{equation}
where $j\ge \al\ge \be$ is such that $t''\in \fin_\al$. Notice that since $j\notin I_1$, and $j\ge
\al\ge\be=\max I_1\cap [0,\dots,j]$, the equation
\begin{equation}
\text{$x_0+T^{k-r}x_1+T^{k-\al}x_2\rel_1' x_0+T^{k-\al}x_2$ is true,}
\end{equation}
for all $r\le j$. Hence,
\begin{equation}
\text{$x_0+T^{k-r}x_1+T^{k-\al}x_2\rel x_0+T^{k-\al}x_2$ is true.}
\end{equation}
There are two now two subcases to consider:

\noindent (a.1) $l\le j $. Then
\begin{equation}
\text{$s'+T^{k-\al}x_2 \rel s'+T^{k-l}x_1+T^{k-\al}x_2\rel t'+T^{k-l'}x_1+T^{k-\al}x_2$ is true,}
\end{equation}
and hence,
\begin{equation}
\text{$x_0+T^{k-l'}x_1+T^{k-l'}x_2+T^{k-\al}x_3\rel x_0+T^{k-l'}x_1+T^{k-\al}x_3$ is true,}
\end{equation}
which implies that
\begin{equation}
\text{$x_0+T^{k-(j+1)}x_1+T^{k-\al}x_2\rel x_0+T^{k-\al}x_2$ holds.}
\end{equation}
By Proposition \ref{red88},
\begin{equation}
\text{ $x_0+T^{k-(j+1)}x_1+T^{k-\al}x_2\rel_1' x_0+T^{k-\al}x_2$  holds,}
\end{equation}
which  contradicts  the fact that $j+1\in I_1$.

\noindent (a.2) $j+1\le l<l'$. Then, the equation
\begin{equation}
\text{$s'+T^{k-l}(x_0+T^{l-j}x_1)+T^{k-\al}x_2
 \rel s'+T^{k-l}x_0+T^{k-\al}x_2$ holds,}
\end{equation}
and hence,
\begin{equation}
\text{$t'+T^{k-l'}x_0+T^{k-(j+l'-l)}x_1+T^{k-\al}x_2
 \rel  t'+T^{k-l'}x_0+T^{k-\al}x_2 $ holds,}
\end{equation}
which implies that
\begin{equation}\label{mmmds}
\text{$x_0+T^{k-(j+l'-l)}x_1+T^{k-\al}x_2\rel  x_0+T^{k-\al}x_2$ holds.}
\end{equation}
Since  $i'-i>0$ the assertion (\ref{mmmds}) contradicts  the fact that $j+1\in I_1$.

\noindent$(b)$ {$j\in I_1$}. We assume that $j+1\in I_1$ because otherwise there is nothing to
prove. Then
$$\rel_1(j+1)=\rel_1(j)\cap \theta^{(1)}_{j+1,l}\cap \max_{j+1},$$ where $l=l^{(1)}_{j+1}$. Suppose that $s\rel
t$. By the inductive hypothesis, $s\rel_1'(j) t$, and in particular $\max_j(s)=\max_j(t)$. Let
$m_0=\max\{\maxi_{j+1}s,\maxi_{j+1}t\}$.  First we show that
\begin{equation}\label{ecua1}
(s\rest[m_0,\maxi_j(s)])^{-1}(l)=
 (t\rest[m_0,\maxi_j(s)])^{-1}(l),
\end{equation}
i.e., for all $n\in [m_0,\max_j(s)]$, $s(n)=l$ iff $t(n)=l$. Suppose  not,  and let
\begin{equation*}
 m_1=\max\{  m\in [m_0,\maxi_j(s)]\, : \, (s(m)=l \text{ or }
t(m)=l ) \text{ and }s(m)\neq t(m)\}.
\end{equation*}
Suppose that $s(m_1)=l$, and that $t(m_1)\neq 0$. Let $n_1$ be the unique integer $n$ such that
 $T^{k-C_{\mathbb{B}}(n)}a_n(m_1)=s(m_1)=l$, and let
$h=C_{\mathbb{B}}(n_1)\ge l$. So, $h'=C_\mathbb{B}(n_1)\neq h$,  $s=s'+T^{k-h}a_{n_1}+s''$,  and
$t=t'+T^{k-h'}a_{n_1}+t''$, with $s'', t''$ both $j$-vectors. By definition of $m_1$, the equation
\begin{equation}
\text{ $x+s''\rel_1'(j+1)x+t''$ holds,}
\end{equation}
and hence,
\begin{equation}
\text{$x+s''\rel_1'x+t''$ and $x+s''\rel x+t''$ also both hold.}
\end{equation}
So, $s'+T^{k-h}a_{n_1}+s''\rel t''+T^{k-h'}a_{n_1}+ s''$, and hence, the equation
\begin{equation}
\text{$s'+T^{k-h}x_0+T^{k-j}x_1\rel t'+T^{k-h'}x_0+T^{k-j}x_1$ holds.}
\end{equation}
There are two subcases to consider:

\noindent (b.1) $h>h'$. Since $x_0+T^{k-r}x_1+T^{k-j}x_2\rel_1' x_0+T^{k-j}x_2$ is true, the
equation $x_0+T^{k-r}x_1+T^{k-j}x_2\rel x_0+T^{k-j}x_2 $ holds  for every $r<l$. Since $l+h'-h<l$,
\begin{align}
s'+T^{k-h}x_0+T^{k-l}x_1+T^{k-j}x_2\sim s'+T^{k-h}(x_0+T^{h-l}x_1)+T^{k-j}x_2\rel \\
t'+T^{k-h'}(x_0+T^{h-l}x_1)+T^{k-j}x_2\sim t'+T^{k-h'}x_0+T^{k-(l+h'-h)}x_1+T^{k-j}x_2\rel \\
\sim t'+T^{k-h'}x_0+T^{k-j}x_2 \rel s'+T^{k-h}x_0+T^{k-j}x_2 \text{ hold.}
\end{align}
Notice that  we have used that $h\ge l$, and so $T^{h-l}$ makes sense. Summarizing, the equation
\begin{equation}
\text{$s'+T^{k-h}x_0+T^{k-l}x_1+T^{k-j}x_2 \rel s'+T^{k-h}x_0+T^{k-j}x_2$ holds,}
\end{equation}
and hence,  the equation
\begin{equation}
\text{$x_0+T^{k-l}x_1+T^{k-j}x_2 \rel_1' x_0+T^{k-j}x_2$ holds,}
\end{equation}
which  is a contradiction with the fact that $\rel_1'\con \theta^{1}_{j+1,l}$.

\noindent (b.2) $h<h'$. Then $h'>l$, and repeating the previous argument used for the case $h>h'$,
we conclude that the equation
\begin{equation}
\text{$t'+T^{k-h'}x_0+T^{k-l}x_1+T^{k-j}x_2 \rel t'+T^{k-h'}x_0+T^{k-j}x_2$ holds,}
\end{equation}
and hence,
\begin{equation}
\text{$x_0+T^{k-l}x_1+T^{k-j}x_2 \rel_1' x_0+T^{k-j}x_2$ holds,}
\end{equation}
which is a contradiction.

The proof will be finished once  we show that $\max_{j+1}s=\max_{j+1}t$. So suppose otherwise,
without loss of generality, that $\max_{j+1}s>\max_{j+1}t$. Let $n_1\in \N$ be such that
$\maxi_{j-1}(s)=\maxi_{j-1}(T^{k-h}b_{n_1})$, where  $h=C_{\mathbb{B}}(n_1)\ge j+1$. Then one has
the decomposition  $s=s'+T^{k-h}a_{n_1}+s''$, $t=t'+T^{k-h'}a_{n_1}+t''$, where $h'<h$ and
$s'',t''$ are $j$-vectors. From (\ref{ecua1}), it follows that
\begin{equation}
\text{$x_0+s''\rel_1'(j+1) x_0+t''$ holds,}
\end{equation}
and hence,
\begin{equation}
s'+T^{k-h}a_{n_1}+t''\rel t'+T^{k-h'}a_{n_1}+t''.
\end{equation}
This implies that the equation
\begin{equation}
\text{$s'+T^{k-h}x_0+T^{k-j}x_1\rel t'+T^{k-h'}x_0+T^{k-j}x_1$ is true.}
\end{equation}
Using a similar argument to the above, we  arrive at  the equation
\begin{equation}
\text{$s'+T^{k-h}(x_0+T^{h-l}x_1)+T^{k-j}x_2\rel t'+T^{k-h'}(x_0+T^{h-l}x_1)+T^{k-j}x_2$  is true,}
\end{equation}
and hence,
\begin{align}
s'+T^{k-h}x_0+T^{k-l}x_1+T^{k-j}x_2\rel
 t'+T^{k-h'}x_0+T^{k-(l+h'-h)}x_1+T^{k-j}x_2\rel \\
\rel t'+T^{k-h'}x_0+T^{k-j}x_2 \rel s'+T^{k-h}x_0+T^{k-j}x_2 \text{ holds,}
\end{align}
which is again a contradiction, since it implies that
\begin{equation}
\text{$x_0+T^{k-l}x_1+T^{k-j}x_2 \rel_1' x_0+T^{k-j}x_2$ holds.}
\end{equation}
\fprue

\prop\label{red7} Suppose that $a,b$ are $k$-vectors of $\mathbb{B}$, $s,t$ are $(\le\! (k-1))$-vectors of
$\mathbb{B}$ such that $a<s$ and $b<t$, and suppose that $a+s\rel b+t$.
\begin{enumerate}
\item If $a<t$ and $b<s$, then $a+s\rel a+t$ and hence $a+t\rel b+t$.
\item If  $a<t$ and $\max_k a<\max_k b$, then $a+s\rel a+t$ and hence $a+t\rel b+t$.
\end{enumerate}
\fprop
\prue
(i)  is a consequence of  Proposition \ref{red3}(1) and Lemma \ref{rel1}. Let us prove (ii). To do
this, suppose that $a,b,s,t$ are as in the statement. By Lemma \ref{rel1} one has that $a+s\sim_1'
b+t$. Since $\max_k(a+s)<\max_k(b+t)$, we have that $\sim_1'=\sim_1 $, where $\sim_1$ is a
$\max$-relation of $\fin_{k-1}$.  This implies that $s\sim_1 t$, from which the desired result
easily follows.
\fprue

\subsection{Determining the relation $\rel$}
We already know that $\rel\con\rel_1'$. The following  identifies  the staircase equivalence
relation that will be equal to $\rel$ on $\mathbb{B}$ in terms of which equations hold or not in
$\mathbb{B}$. This will conclude the proof of Theorem \ref{maint}.
\teor \label{main11}
\begin{enumerate}
\item[]
\item Suppose that $x_0+T^{k-(l-1)} x_1 +x_2\rel x_0+ x_2$ is true, and
$x_0+T^{k-l} x_1 +x_2\rel x_0+ x_2$ is false.
 Then $\rel=\rel_0'\cap \rel_{\theta_{l}^2}\cap \rel_1'$.
 \item Suppose that $x_0+x_1+x_2\rel x_0+x_2$ is true.
\begin{enumerate}
\item If $Tx_0+x_1+x_2\rel Tx_0+x_2$ and $x_0+x_1+Tx_2\rel x_0+Tx_2$ are both false, then
$\rel=\rel_0'\cap \min_k\cap \max_k\cap \rel_1'$.
\item If $Tx_0+x_1+x_2\rel Tx_0+x_2$ is true, and  $x_0+x_1+Tx_2\rel x_0+Tx_2$ is false,
then $\rel=\rel_0'\cap \max_k\cap \rel_1'$.
\item If $Tx_0+x_1+x_2\rel Tx_0+x_2$ is false, and  $x_0+x_1+Tx_2\rel x_0+Tx_2$ is true,
then $\rel=\rel_0'\cap \min_k\cap \rel_1'$.
\item If $Tx_0+x_1+x_2\rel Tx_0+x_2$ and  $x_0+x_1+Tx_2\rel x_0+Tx_2$ are both  true,
then $\rel=\rel_0'\cap \rel_1'$.

\end{enumerate}

\end{enumerate}

\fteor
The proof is done in various steps. (i) is  in Corollary \ref{case1a}, and (ii.a), (ii.b), (ii.c)
and (ii.d) in Corollary \ref{case2a}, and  Lemmas \ref{case2b}, \ref{case2c} and \ref{case2d}
respectively. We start with the following proposition that gives one of the inclusions.

\prop\label{red10}
\begin{enumerate}
\item[] \item If the equation $Tx_0+x_1+x_2\rel Tx_0+x_2$ is true,
then $\rel_0'\cap \max_k \cap \rel_1'\con \rel $.

\item If the equation $x_0+x_1+Tx_2\rel x_0+Tx_2$ is true, then
$\rel_0'\cap \min_k \cap \rel_1'\con \rel $. \item If the equation $x_0+T^{k-(l-1)} x_1 +x_2\rel
x_0+ x_2$ is true, then $\rel_0'\cap \rel_{\theta^2_l}\cap \rel_1'\con \rel $, for every $l\le k$.

\item If the equation $x_0+x_1+x_2\rel x_0+x_2$ is true, then
$\rel_0'\cap \min_k \cap \max_k\cap \rel_1'\con \rel$. \item If the equations $Tx_0+x_1+x_2\rel
Tx_0+x_2$ and  $x_0+x_1+Tx_2\rel x_0+Tx_2$ are both  true, then $\rel_0'\cap \rel_1'\con \rel$.
\end{enumerate}
\fprop
\prue
\noindent (i): Suppose that the equation $Tx_0+x_1+x_2\rel Tx_0+x_2$ holds. Then,
$Tx_0+Tx_1+x_2+x_3\rel Tx_0+x_3$ holds, and the equations $Tx_0+Tx_1+x_2+x_3\rel Tx_0+x_3\rel
Tx_0+x_2+x_3$ also hold. This implies that the equation
\begin{equation}
\text{$Tx_0+Tx_1+x_2\rel Tx_0+x_2$ is true.}
\end{equation}
Hence the relation $\rel_0$ is a $\min$-relation, which implies that $\rel_0'$ so is a
$\min$-relation. Set $R=\rel_0'\cap \max_k \cap \rel_1'$ and suppose that $sRt$. Then $\max_k s
=\max_k t$. Let $n$ be  such that $\max_k s=\max_k=\max_k b_n$. Therefore, $s=s'+a_{3n+1}+s''$ and
$t=t'+a_{3n+1}+t''$. It is not difficult to show that the equation
\begin{equation}
\text{$Tx_0+x_1+x_2 \, R \, Tx_0+x_2$ holds.}
\end{equation}
So, we may assume  that $s'$ and $t'$ are $(k-1)$-vectors of $\mathbb{A}$. Since $s \rel_1' t$, we
have that $s'+a_{3n+1}+s''\rel s'+a_{3n+1}+t''$. Since $s \rel_0' t$, we have that
$s'+a_{3n+1}+t'' \rel_0' t'+a_{3n+1}+t''$, and hence $s'\rel_0 t'$, which implies that $s'+x \rel
t' +x$ is true. In particular $s'+a_{3n+1}+t'' \rel t'+a_{3n+1}+t''$, i.e., $s\rel t$.

The proofs of (ii), (iii) and (iv) are  similar. We leave the details to the reader. Let us check
point (v): Fix $s=f_\mathbb{A}s+a_{n_0}+m_\mathbb{A}s+a_{n_1}+l_{\mathbb{A}}s$,
$t=f_\mathbb{A}t+a_{m_0}+m_\mathbb{A}t+a_{m_1}+l_{\mathbb{A}}t$ such that $s R t$, where
$R=\rel_0'\cap \rel_1'$.  If $m_0=n_0$, then $sR\cap\min_k t$, and hence we are done by 2. So,
suppose that $n_0<m_0$. Since $\rel'_0$ is a $\min$-relation and $\rel_1'$ is a $\max$-relation,
the equations $Tx_0+x_1+x_2 R Tx_0+x_2$ and  $x_0+x_1+Tx_2 R x_0+Tx_2$ are true. Therefore, $s R
f_\mathbb{A}s+a_{n_0}+l_{\mathbb{A}}s$ and $t R f_\mathbb{A}t+a_{m_0}+l_{\mathbb{A}}t$. Since
$s\rel f_\mathbb{A}s+a_{n_0}+l_{\mathbb{A}}s$ and $t \rel f_\mathbb{A}t+a_{m_0}+l_{\mathbb{A}}t$
the proof will be finished if we show that
\begin{equation}
f_\mathbb{A}s+a_{n_0}+l_{\mathbb{A}}s\rel f_\mathbb{A}t+a_{m_0}+l_{\mathbb{A}}t.
\end{equation}  Since $f_\mathbb{A}s+a_{n_0}+l_{\mathbb{A}}s\rel_0' f_\mathbb{A}t+a_{m_0}+l_{\mathbb{A}}t$ and
$f_\mathbb{A}s+a_{n_0}+l_{\mathbb{A}}s\rel_1' f_\mathbb{A}t+a_{m_0}+l_{\mathbb{A}}t$,  by the last
point of Proposition \ref{red3} (for both $\rel_0'$ and $\rel_1'$), we have that
 \begin{equation}
\text{$f_\mathbb{A}s+a_{m_0}+l_{\mathbb{A}}t\rel f_\mathbb{A}t+a_{m_0}+l_{\mathbb{A}}t$ and
$f_\mathbb{A}s+a_{n_0}+l_{\mathbb{A}}s\rel f_\mathbb{A}s+a_{n_0}+l_{\mathbb{A}}t$.}
\end{equation}
 But
$f_\mathbb{A}s+a_{m_0}+l_{\mathbb{A}}t\rel f_\mathbb{A}s+a_{n_0}+l_{\mathbb{A}}t$, and we are done.
\fprue

\lema\label{max_k}
If the equation $x_0+x_1+T x_2\rel x_0+Tx_2$  is false, then  $\rel\con \maxi_k$.
\flema
\prue
Suppose that $s\rel t$ but $\maxi_{k}s>\maxi_k t$. Set
\begin{align*}
s=& f_\mathbb{A}s+a_{n_0}+m_\mathbb{A}s+a_{n_1}+l_\mathbb{A}s \\
t=& f_\mathbb{A}t+a_{m_0}+m_\mathbb{A}t+a_{m_1}+l_\mathbb{A}t,
\end{align*}
where $n_1>m_1$. Set $l_\mathbb{A}t=t'+T^{k-i}a_{n_1}+t''$, where $t'<T^{k-i}a_{n_1}<t''$, and
$i<k$. By Proposition \ref{red7}
\begin{equation}
f_\mathbb{A}t+a_{m_0}+m_\mathbb{A}t+a_{m_1}+t'+T^{k-i}a_{n_1}+l_\mathbb{A}s\rel
f_\mathbb{A}s+a_{n_0}+m_\mathbb{A}s+a_{n_1}+l_\mathbb{A}s,
\end{equation}
and therefore, the equation
\begin{equation}
\text{$f_\mathbb{A}t+a_{m_0}+m_\mathbb{A}t+a_{m_1}+t'+T^{k-i}x_0+Tx_1\rel
f_\mathbb{A}s+a_{n_0}+m_\mathbb{A}s+x_0+T x_1$ holds.}
\end{equation}
Since $\sim\con \sim_1'$ and $\sim_1'$ is a canonical relation, the $\sim_1'$-equation
\begin{equation}
\text{$f_\mathbb{A}t+a_{m_0}+m_\mathbb{A}t+a_{m_1}+t'+T^{k-i}x_0+Tx_1\rel_1'
f_\mathbb{A}s+a_{n_0}+m_\mathbb{A}s+x_0+T x_1$ holds.}\end{equation} Since $\rel_1'$ is a
staircase relation, the truth of the last equation  implies that $k\notin I_1(\rel_1')$, and hence
$\rel_1'$ is a $\max$-relation with $\max (I_1(\rel_1'))$ at most $k-1$. Therefore,
\begin{equation}
\text{$f_\mathbb{A}t+a_{m_0}+m_\mathbb{A}t+a_{m_1}+t'+T^{k-i}x_0+Tx_1\rel_1'
f_\mathbb{A}t+a_{m_0}+m_\mathbb{A}t+a_{m_1}+t'+Tx_1$ is true,}
\end{equation}
which implies that
\begin{equation}
\text{$f_\mathbb{A}t+a_{m_0}+m_\mathbb{A}t+a_{m_1}+t'+T^{k-i}x_0+Tx_1\rel
f_\mathbb{A}t+a_{m_0}+m_\mathbb{A}t+a_{m_1}+t'+Tx_1$ is true.}
\end{equation}
Hence, the equation
\begin{equation}
\text{$f_\mathbb{A}t+a_{m_0}+m_\mathbb{A}t+a_{m_1}+t'+Tx_1\rel_1'
f_\mathbb{A}s+a_{n_0}+m_\mathbb{A}s+x_0+T x_1$ holds,}
\end{equation}
from which we conclude that
\begin{equation}
\text{$x_0+x_1+Tx_2\rel x_0+Tx_2$ is true,}
\end{equation}
a contradiction.
\fprue
\lema\label{med_i} Suppose that $x_0+T^{k-(l-1)} x_1 +x_2\rel x_0+ x_2$ is true but $x_0+T^{k-l}
x_1+x_2\rel x_0+ x_2$ is false. Then $\rel\con \rel_{\theta^{2}_{l}}$. In particular, $\rel\con
\mini_k \cap \maxi_k$.
\flema
\prue Fix  $l$ as in the statement. Since we assume that the equation
\begin{equation}
\text{$x_0+T^{k-l}x_1+x_2\sim x_0+x_2$ is false,}
\end{equation}
by Proposition \ref{preq}(1,2), we know that
\begin{equation}
\text{ $x_0+x_1+Tx_2\rel x_0+Tx_2$ is false.}
\end{equation}
So, by Lemma \ref{max_k}, we obtain that $\rel\con \max_k$. Suppose that $s\sim t$. Take the
decomposition
\begin{align*}
s=& f_\B s + b_{n_0}+m_\B s+ b_{m} +l_\B s \\
t=& f_\B t + b_{n_1}+m_\B t+ b_{m} +l_\B s,
\end{align*}
where we implicitly assume that $l_\B s=l_\B t$, since $s\sim_1't$. Observe that  showing that
$s\theta_l^2 t$ is the same that proving that
\begin{equation}\label{asdfdw}
\text{ for all $n\in [\min\{n_0,n_1\},m]$, either $C_\B(s)(n), C_\B(t)(n)<l$, or
$C_\B(s)(n)=C_\B(t)(n)$.}
\end{equation}
Assume on the contrary that (\ref{asdfdw}) is false, and let $\al$ be the last $n\in
[\min\{n_0,n_1\},m]$ for which
\begin{equation}
\text{  $\max\{C_\B(s)(n),C_\B(t)(n)\}\ge l$ and $C_\B(s)(n)\neq C_\B(t)(n)$.}
\end{equation}
Set $l_0=C_\B(s)(\al)$, and $l_1=C_\B(s)(\al)$. Notice that $\al<m$. Without loss of generality,
we assume that $l_1<l_0$ (the other case has a similar proof). Set
\begin{align*}
s'=&\sum_{n<\al}T^{k-C_\B(s)(n)}b_n\\
t'=&\sum_{n<\al}T^{k-C_\B(t)(n)}b_n.
\end{align*}
Using this notation, we have that the equation
\begin{equation}\label{iueww}
\text{ $s'+T^{k-l_0}x_0+x_1\sim t'+T^{k-l_1}x_0+x_1$ holds.}
\end{equation}
 There are two cases:

\noindent \underline{$n_0\le n_1$}. We first show that in this case  $s'+T^{k-l_0}x_0$ is a
$k$-term. If $n_0=n_1$, then $\al>n_0$, and hence $s'$ is a $k$-vector. Suppose that $n_0<n_1$. If
$\al>n_0$, then $s'$ is a $k$-term. If $\al=n_0$, then $l_0=k$, and clearly $s'+T^{k-k}x_0=s'+x_0$
is a $k$-term. We consider two subcases:

\noindent \emph{(a)} $l_1<l \le l_0$. Then, by our assumption that $x_0+T^{k-(l-1)}x_1+x_2\sim
x_0+x_2$ holds, we have that
\begin{equation}
\text{$s'+T^{k-l_0}x_0+T^{k-l_1}x_1+ x_2 \sim  s'+T^{k-l_0}x_0+ x_2$ holds.}
\end{equation}
By (\ref{iueww}),
\begin{align}
& s'+T^{k-l_0}x_0+T^{k-l_1}x_1+ x_2 \sim  t'+T^{k-l_1}x_0+T^{k-l_1}x_1+ x_2\sim \\
& \sim s'+T^{k-l_0}x_0+T^{k-l_0}x_1+ x_2\text{ holds,}
\end{align}
which implies that the equation
\begin{equation}
\text{$s'+T^{k-l_0}x_0+T^{k-l_0}x_1+ x_2 \sim  s'+T^{k-l_0}x_0+ x_2$ holds.}
\end{equation}
This contradicts  the fact that $l_0\ge l$.

\noindent \emph{(b)} $l\le  l_1 < l_0$. Then,
\begin{equation}
\text{$s'+T^{k-l_0}x_0+ T^{k-l_1}(T^{l_0-l})x_1 +x_2\sim s'+T^{k-l_0}x_0+x_2$ holds,}
\end{equation}
and by (\ref{iueww}),
\begin{align}
& s'+T^{k-l_0}x_0+ T^{k-l_1}(T^{l_0-l})x_1 +x_2\sim  t'+T^{k-l_1}x_0+ T^{k-{l_1}}(T^{l_0-l})x_1 +x_2 \sim
\\
& \sim s'+T^{k-l_0}x_0+T^{k-l}x_1+x_2
 \text{ holds.}
\end{align}
Again, this yields  a contradiction.

\noindent \underline{$n_1<n_0$}. It can be shown that $t'+T^{k-l_1}x_0$ is a $k$-term.  We
consider the same two subcases as above:

\noindent \emph{(a)} $l_1<l \le l_0$. Then
\begin{equation}
\text{$t'+T^{k-l_1}x_0+T^{k-l_1}x_1+ x_2 \sim  s'+T^{k-l_1}x_0+ x_2$ holds,}
\end{equation}
and hence,
\begin{equation}
\text{$s'+T^{k-l_0}x_0+T^{k-l_0}x_1+ x_2 \sim  s'+T^{k-l_0}x_0+ x_2$ holds,}
\end{equation}
which, by (\ref{iueww}), implies that
\begin{equation}
\text{$t'+T^{k-l_1}x_0+T^{k-l_0}x_1+ x_2 \sim  t'+T^{k-l_1}x_0+ x_2$ holds,}
\end{equation}
a contradiction, since $l_0\ge l$.

\noindent \emph{(b)} $l\le  l_1 < l_0$. Then
\begin{equation}
\text{$t'+T^{k-l_1}x_0+ T^{k-l_1}(T^{l_0-l})x_1 +x_2\sim t'+T^{k-l_1}x_0+x_2$ holds.}
\end{equation}
Using that
\begin{align}
& t'+T^{k-l_1}x_0+ T^{k-l_1}(T^{l_0-l})x_1 +x_2\sim s'+T^{k-l_0}x_0+ T^{k-l}x_1+ x_2 \sim \\
& \sim t'+T^{k-l_1}x_0 +T^{k-l}x_1+ x_2 \text{ holds,}
\end{align}
we arrive at a contradiction.
\fprue

\cor\label{case1a}
 Suppose that $x_0+T^{k-(l-1)} x_1 +x_2\rel x_0+ x_2$ is true, but $x_0+T^{k-l}
x_1+x_2\rel x_0+ x_2$ is false. Then, $\rel=\rel_0'\cap\rel_{\theta^2_l}\cap \rel_1'$. \fcor
\prue
By Proposition \ref{red10}, $\rel_0'\cap\rel_{\theta^2_l}\cap \rel_1'\con \rel$.  We only need to
show that $\rel\con \rel_0'$. Suppose that $s\rel t$, and  consider the decomposition
\begin{align*}
s=& f_\A s+a_{n_0}+m_\A s+a_{m_0}+l_\A s \\
t=& f_\A t+a_{n_1}+m_\A t+a_{m_1}+l_\A t.
\end{align*}
Since $\max_k s=\max_k t$, we have that $m_0=m_1$, and since $s\sim_1't $, by Proposition
\ref{stairr2}(4), we may assume that $l_\A s\sim_1 l_\A t$. By Lemma  \ref{med_i},
$s\rel_{\theta_2^l} t$, and using the fact that the equations $x_0+T^{k-j} x_1 +x_2\rel x_0+ x_2$
are true for all $j<l$, we may also assume that $n_0=n_1$ and $m_\A s=m_\A t$. Therefore, the
equation $f_\A s+x_0\rel f_\A t +x_0$ holds. By definition of $\sim_0$, we have that $f_\A s\sim_0
f_\A t$, and by Proposition \ref{stairr2}(3), $s\sim_0' t$, as desired.
\fprue

\lema\label{case2b}
Suppose that $Tx_0+x_1+x_2\rel Tx_0+x_2$ is true, and $x_0+x_1+Tx_2\rel x_0+Tx_2$ is false. Then,
$\rel=\rel_0'\cap \max_k\cap \rel_1'$.
\flema
\prue
We only need to  show that $\rel\con \rel_0'$. Suppose that $s\rel t$. Consider the following decompositions
of $s$ and $t$
\begin{align*}
s= & f_\A s+a_{n_0}+m_\A s+a_{m}+l_\A s \\
t= & f_\A t+a_{n_1}+m_\A t+a_{m}+l_\A s.
\end{align*}
Notice that,  since $Tx_0+x_1+x_2\rel Tx_0+x_2$ is true, we have that $x_0+x_1+x_2\rel x_0+x_2$ is
true. Hence, we may assume that $m_\A s=m_\A t=0$. Notice also that, since
\begin{equation}
\text{$Tx_0+x_1+x_2\rel Tx_0+x_2$ is true,}
\end{equation}
and since $f_\A s$ and $f_\A t$ are $(k-1)$-vectors (this is why we use the decompositions of
vectors of $\B$ in $\A$), we have that
\begin{equation}
\text{$s\rel f_\A s +a_{n_2}+l_\A s$ and   $t\rel f_\A t +a_{n_2}+l_\A s$.}
\end{equation}
This implies that $f_\A s\sim_0 f_\A t$, and, by Proposition \ref{stairr2}(1,3),
\begin{equation}
s\sim_0' f_\A s +a_{n_2}+l_\A s \sim_0' f_\A t +a_{n_2}+l_\A t\sim_0' t,
\end{equation}
as desired.\fprue

\prop\label{nomax_k} Suppose that $x_0+x_1+Tx_2\rel x_0+Tx_2$ holds, and suppose that
$Tx_0+x_1+x_2\rel Tx_0+x_2$ is false. Then, $\rel\con \min_k$.
\fprop
\prue
Suppose that $s\rel t$. Take the decomposition
\begin{align*}
s= & f_\A s+a_{n_0}+m_\A s+a_{m_0}+l_\A s \\
t= & f_\A t+a_{n_1}+m_\A t+a_{m_1}+l_\A t.
\end{align*}
Suppose that $n_0\neq n_1$, and without loss of generality assume that $n_0<n_1$.  Since
$x_0+x_1+Tx_2\rel x_0+Tx_2$ holds, we have that
\begin{equation}
f_\A s+a_{n_0}+l_\A s\rel f_\A t+a_{n_1}+l_\A t.
\end{equation}
By Proposition \ref{red7}(2), we have that
\begin{equation}
f_\A s+a_{n_0}+l_\A t\rel f_\A t+a_{n_1}+l_\A t,
\end{equation}
and hence (since $l_\A t$ is a $(k-1)$-vector), the equation
\begin{equation}
\text{$f_\A s+x_0+Tx_2\rel f_\A t+x_1+Tx_2 $ holds.}
\end{equation}
This implies that
\begin{equation}
\text{$f_\A s+x_0+x_1+Tx_3\rel f_\A t+x_2+Tx_3\rel f_\A s+x_1+Tx_3 $ holds,}
\end{equation}
which implies that the equation
\begin{equation}
\text{$f_\A s +x_0+x_2\rel f_\A s +x_2$ is true.}
\end{equation}
Since $f_\A s$ is a $(k-1)$-vector, we have that
\begin{equation}
\text{$T x_0 +x_1+x_2\rel T x_0 +x_2$ is true,}
\end{equation}
a contradiction.
\fprue

\lema\label{case2c}
Suppose that $x_0+x_1+Tx_2\rel x_0+Tx_2$ is true, and $Tx_0+x_1+x_2\rel Tx_0+x_2$ is false. Then,
$\rel=\rel_0'\cap \min_k\cap \rel_1'$.
\flema
\prue
By Proposition \ref{red10}, we have that $\rel_0'\cap \min_k\cap \rel_1'\con \rel$.  Let us show
 that $\rel\con\rel_0'\cap \min_k\cap \rel_1'$. By Proposition \ref{nomax_k} and Lemma
 \ref{rel1},
 we have that  $\rel \con \min_k\cap \rel_1'$. So, we only need to show that $\rel\con \rel_0'$.
 Suppose that $s\rel t$ with
\begin{align*}
s= & f_\A s+a_{n_0}+m_\A s+a_{n_1}+l_\A s \\
t=& f_\A t+a_{n_0}+m_\A t+a_{m_1}+l_\A t.
\end{align*}
Since the equation $x_0+x_1+Tx_2\rel x_0+Tx_2$ is true, we have that
\begin{equation}
f_\A s +a_{n_0}+l_\A s\sim f_\A t +a_{n_0}+l_\A t,
\end{equation}
and, by Proposition \ref{red7},
\begin{equation}
f_\A s +a_{n_0}+l_\A s\sim f_\A t +a_{n_0}+l_\A s,
\end{equation}
which easily leads to that $s\rel_0' t$.
\fprue

\lema\label{nomed_i}
Suppose that $x_0+x_1+x_2\sim x_0+x_2$ is true , and that $x_0+x_1+Tx_2\rel x_0+Tx_2$ and
$Tx_0+x_1+x_2\rel Tx_0+x_2$ are both false. Then, $\rel\con \min_k\cap \max_k$.
\flema
\prue
By Lemma \ref{max_k}, we know that $\rel\con \max_k$, and by Lemma  \ref{rel1}, $\rel\con \rel_1'$.
So, we only need to show that $\sim \con \min_k$. Suppose that $s\rel t$, set
\begin{align*}
s=& f_\A s+a_{n_0}+m_\A s+a_{m}+l_\A s \\
t= & f_\A t+a_{n_1}+m_\A t+a_{m}+l_\A t.
\end{align*}
Suppose on the contrary that $n_0<n_1$. There are two cases to consider:

\noindent \underline{$n_1=m$}. Hence, $n_0<m $ and
\begin{equation}
s\rel f_\A s +a_{n_0}+a_{m}+l_\A s \text{ and } t=f_\A t+a_{m}+l_\A t.
\end{equation}
By Proposition \ref{red7},
\begin{equation}
f_\A s +a_{n_0}+a_{m}+l_\A s \rel f_\A t + a_m +l_\A s,
\end{equation}
which implies that the equation
\begin{equation}
\text{$f_\A s+x_0+x_1\rel f_\A t+x_1 $ is true,}
\end{equation}
a contradiction, since  $f_\A s$ is a $(k-1)$-vector.

\noindent \underline{$n_1<m$}. Then, by our assumptions, and Proposition \ref{red7},
\begin{equation}
f_\A s+a_{n_0}+a_{m}+l_\A s \rel f_\A t+a_{n_1}+a_{m}+l_\A s.
\end{equation}
Hence, the equation
\begin{equation}
\text{$f_\A s+x_0+x_2 \rel f_\A t+x_1+x_2$ is true,}
\end{equation}
which readily implies that $Tx_0+x_1+x_2\rel Tx_0+x_2$ must be true, a contradiction.
 \fprue
\cor\label{case2a}
Suppose that $x_1+x_2+x_3\sim x_1+x_3$ is true , and that $x_1+x_2+Tx_3\rel x_1+Tx_3$ and
$Tx_1+x_2+x_3\rel Tx_1+x_3$ are both false. Then, $\rel= \rel_0'\cap\min_k\cap \max_k\cap
\rel_1'$. \fcor
\prue
By Proposition \ref{red10}, $\rel_0'\cap\min_k\cap \max_k\cap \rel_1'\con \rel$. Let us show the
opposite inclusion. By  Lemma  \ref{nomed_i}, we have that $\rel\con \min_k\cap \max_k$. It remains
to show that $\rel\con \rel_0'$. Suppose that $s\rel t$, where $s=f_\A s+a_{n}+m_\A s+a_{m}+l_\A s$
and $t=f_\A t+a_{n}+m_\A t+a_{m}+l_\A s$ (we may assume that $l_\A s=l_\A t$, since $\max_k
(s)=\max_k(t)$). There are two cases: \underline{$n<m$}. Then, $f_\A s+a_{n}+a_{m}+l_\A s \rel f_\A
t+a_{n}+m_\A t+a_{m}+l_\A s$ which directly implies that $s\rel_0' t$.  The proof for
\underline{$n_0=m$} is quite similar.
\fprue

\lema
\label{case2d} Suppose that $Tx_0+x_1+x_2\rel Tx_0+x_2$ and  $x_0+x_1+Tx_2\rel x_0+Tx_2$ are both
true. Then, $\rel=\rel_0'\cap \rel_1'$.
\flema
\prue
It is enough to show that $\rel\con \rel_0'$. Suppose that $s\rel t$, with
$s=f_\mathbb{A}s+a_{n_0}+m_\mathbb{A}s+a_{m_0}+l_{\mathbb{A}}s$ and
$t=f_\mathbb{A}t+a_{n_1}+m_\mathbb{A}t+a_{m_1}+l_{\mathbb{A}}t$. We may assume that
$s=f_\mathbb{A}s+a_{n_0}+l_{\mathbb{A}}s$, and $t=f_\mathbb{A}t+a_{n_1}+l_{\mathbb{A}}t$. W.l.o.g.
we assume that $n_0\le n_1$, and hence, by Proposition \ref{red7},
\begin{equation}\label{mklwqq}
f_\mathbb{A}s+a_{n_0}+l_{\mathbb{A}}t \sim f_\mathbb{A}t+a_{n_1}+l_{\mathbb{A}}t.
\end{equation}
\noindent Case {$n_0=n_1$}. By definition of $\rel_0'$, (\ref{mklwqq}) implies that
\begin{equation}
f_{\mathbb{A}}t+a_{n_0}+l_\mathbb{A}t\rel_0' f_\mathbb{A}s+a_{n_0}+l_{\mathbb{A}}t,
\end{equation}  but trivially
$f_\mathbb{A}s+a_{n_0}+l_{\mathbb{A}}t\rel_0' f_\mathbb{A}s+a_{n_0}+l_{\mathbb{A}}s$, and we are
done.

\noindent Case  {$n_0<n_1$}. Then,
\begin{equation}
\text{$f_\A s+x_0+T x_2\sim f_\A t +x_1+Tx_2$ is true,}
\end{equation}
which easily yields
\begin{equation}
\text{$f_\A s+x_1+T x_3\sim f_\A t +x_1+Tx_3$ is true.}
\end{equation}
This implies that $s\sim_0' t$.\fprue

\cor
Every equivalence relation on $\fin_k$ is canonical in some sos. \qed
 \fcor
This corollary has the following local version.
\cor
For every  block sequence $A$ and every equivalence relation $\rel$ on $\langle A \rangle$ there is
an sos $B\in [A]^{[\infty]}$ on which $\rel$ is canonical. \fcor
\prue
Fix the canonical isomorphism $\La:\fin_k\to \langle  A \rangle$ (i.e., the extension of $\Theta
e_n\mapsto a_n$).  It is not difficult to show the following facts:
\begin{enumerate}
\item $B=(b_n)_n$ is an sos iff $FB=(Fb_n)_n$ is an sos.
\item For every canonical equivalence  relation $\rel_{can}$, every sos $B$, and $s,t\in \langle B \rangle$,
$s\rel_{can}t$ iff $F^{-1}s \rel_{can} F^{-1}t$.
\end{enumerate}
We define $\rel'$ on $\fin_k$ by $s\rel' t$ iff $Fs\rel Ft$. Find a canonical equivalence relation
$\rel_{can}$ and an sos $B$ such that $\rel$ and $\rel_{can}$ are the same on $\langle  B \rangle$.
Let $C=FB$, which is an sos. Then $\rel$ and $\rel_{can}$ are the same in $\langle C \rangle$:
$s\rel_{can} t $ iff $F^{-1} s\rel_{can} F^{-1}t$ iff $F^{-1} s\rel' F^{-1}t$ iff $s\rel t$.
\fprue

\cor
Every canonical equivalence relation is a staircase equivalence relation. \fcor
\prue
Notice that, since $\sim$ is canonical in $A$,  $\mathbb{A}=A$ works for both Lemmas \ref{eqnde}
and \ref{reg1}. Hence,  $\rel$ is a  staircase equivalence relation in
$\mathbb{B}=(Ta_{3n}+a_{3n+1}+Ta_{3n+2})_n$. Let $\sim'$ be this staircase relation, which is equal
to $\sim$ when restricted to $\B$. We  show that $\sim$ and $\sim'$ are not only equal in
$\mathbb{B}$, but also in $A$.   Fix $s$ and $t$ in $A$, and take their canonical decompositions in
$A$
\begin{equation*}
\text{ $s=\sum_{n\ge 0}T^{k-C_A(s)(n)}a_n$ and $t=\sum_{n\ge 0}T^{k-C_A(t)(n)}a_n$.}
\end{equation*} Suppose first that $s\rel t $.
Since $\sim$ is canonical, the equation
\begin{equation}
\text{ $\sum_{n\ge 0}T^{k-C_A(s)(n)}x_n\rel \sum_{n\ge 0}T^{k-C_A(t)(n)}x_n$ holds in $A$,}
\end{equation}
and hence, also  in $\mathbb{B}$, i.e.,
\begin{equation}\label{fgtwe}
\text{ $\sum_{n\ge 0}T^{k-C_A(s)(n)}x_n\rel' \sum_{n\ge 0}T^{k-C_A(t)(n)}x_n$ holds in
$\mathbb{B}$.}
\end{equation}
But since $\rel'$ is staircase, it is canonical (Proposition \ref{st->can}), and hence, equation
(\ref{fgtwe}) also  holds in $A$, and in particular, $s\rel' t$.

Suppose now that $s\rel' t$.  Since $\rel'$ is canonical in any sos, the equation
\begin{equation}
\sum_{n\ge 0}T^{k-C_A(s)(n)}x_n\rel' \sum_{n\ge 0}T^{k-C_A(t)(n)}x_n \text{ holds in $A$,}
\end{equation}
hence, also in $\mathbb{B}$. By definition,   $\rel'$ is  equal to $\rel$ restricted to
$\mathbb{B}$, and hence
\begin{equation}\label{loiuty}
\sum_{n\ge 0}T^{k-C_A(s)(n)}x_n\rel \sum_{n\ge 0}T^{k-C_A(t)(n)}x_n \text{ holds in $\mathbb{B}$.}
\end{equation}
Since $\rel$ is canonical, the equation (\ref{loiuty}) holds in $A$, and in particular, $s\sim t$.
\fprue

\section{Counting}\label{counting}
The purpose now is to give an explicit formula for the number $t_k$ of staircase equivalence
relations on $\fin_k$. To do this, recall that  $e_n(1)=\sum_{j=0}^{n}\frac{1}{j!}$ is the
exponential sum-function and that $\Gamma(a,x)=\int_x^\infty t^{a-1} e^{-t}dt$ is the incomplete
Gamma function.  Recall also that $\Gamma(n,1)=(n-1)!e^{-1}e_{n-1}(1)$ for every integer $n$.

Let $\mathcal{A}_k$, $\mathcal{B}_k$ be the set of $\min$-relations and $\max$-relations
respectively, and set $a_k=|\mathcal{A}_k|$ and $b_k=|\mathcal{B}_k|$. Let $\mathcal{C}_k\con
\mathcal{A}_k$ be the set of $\min$-relations $R$ such that $k\notin I_0(R)$, and let
$\mathcal{D}_k\con \mathcal{B}_k$ be the set of $\max$-relations $R$ such that $k\notin I_1(R)$.
Set $c_k=|\mathcal{C}_k|$ and $d_k=|\mathcal{D}_k|$. Notice that
\begin{enumerate}
\item $c_k=a_{k-1}$,
\item $\mathcal{A}_k=\mathcal{A}_{k-1}\cup \conj{R\cap \rel_{\min_k}}{R\in \mathcal{C}_{k-1}}
\cup\{R\cap \rel_{\min_k}\cap \rel_{\theta^0_{k,l}} \, : \, l=-1 \text{ or }l=1,\dots,k-1,\, R\in
\mathcal{A}_{k-1}\setminus \mathcal{C}_{k-1}\}$. So, $a_k=a_{k-1}+c_{k-1}+k(a_{k-1}-c_{k-1})$.

\end{enumerate}
Hence,
\begin{equation}
a_k=(k+1)a_{k-1}-(k-1)a_{k-2}, \,a_0=1, \, a_1=2.
\end{equation}
By standard methods, we conclude that
\begin{equation}\label{kestse}
a_k=\frac{e\,\left( 1 + k \right) \,k!\Gamma(1 + k,1)}{\Gamma(2 + k)}=k! e_k(1).
\end{equation}
Now let $\mathcal{T}_k$ be the set of staircase equivalence relations of $FIN_k$ and
$t_k=|\mathcal{T}_k|$. Then,
\begin{align}
\mathcal{T}_k= &\left(\conj{R\cap S}{R\in \mathcal{A}_k, S\in \mathcal{B}_k}\setminus
\conj{R\cap
S}{R\in \mathcal{A}_k\setminus \mathcal{C}_k, S\in \mathcal{B}_k\setminus \mathcal{D}_k} \right) \\
& \cup\conj{R\cap S\cap \rel_{\theta^2_l}}{R\in \mathcal{A}_k\setminus \mathcal{C}_k, S\in
\mathcal{B}_k\setminus \mathcal{D}_k, l=-1 \text{ or } l=1,\dots,k}.
\end{align}
Hence,
\begin{equation}\label{weweradfd}
t_k= a_k^2-(a_k-c_k)^2+(k+1)(a_k-c_k)^2=k(a_k-a_{k-1})^2+a_{k}^2
\end{equation}
and from (\ref{kestse}) and (\ref{weweradfd}), we obtain that
\begin{equation}
t_k= (k! e_k(1))^2+ k\left(k! e_k(1)-(k-1)!e_{k-1}(1)\right)^2,
\end{equation}
or, equivalently,
\begin{equation}
t_k= e^2 \left[   k\left[\Gamma(k, 1) - \Gamma(k+1, 1) \right]^2 +
        \Gamma(k+1,1)^2\right].
\end{equation}
This is a table with the first few values of $t_k$:
$$\begin{array}{|l|l|l|l|l|l|l|l|}\hline
k & 0 & 1 & 2 &3 & 4 &5 &6 \\
\hline t_k & 1 & 5& 43& 619& 13829& 446881& 19790815\\ \hline
\end{array} $$

\nota Let us say that a canonical equivalence relation $R$ is \emph{linked free} iff $I_0(R)$ and
$I_1(R)$ have no consecutive members and $k\notin I_0(R)\cap I_1(R)$.  The number $l_k$ of linked
free canonical equivalence relations of $\fin_k$ is the Fibonacci number $F_{2k+2}$ for $2k+2$,
since $F_{l+2}$ is the number of subsets of $\{1,2,\dots,l\}$ with no consecutive elements, and
since $R$ is linked free iff  the set $I_0(R)\cup \conj{2k+1-i}{i\in I_1(R)}\con \{1,2,\dots,2k\}$
has no consecutive numbers. \fnota

\section{the finite version}

\teor\label{finite}
 For every $m$ there is some $n=n(m)$ such that  for every equivalence relation $\rel$
 on $\langle e_0,\dots,e_n \rangle$ there is some sos $(a_0,\dots,a_{m-1})\pe( e_0,\dots,e_n)$
  such that $\rel$ is a staircase equivalence relation in $\langle a_0,\dots,a_{m-1} \rangle$.
\fteor
\prue Suppose not. Then, there is some $m$ such that for every $n$ there is  some
equivalence relation $\rel_n$ on $\langle e_0,\dots,e_n \rangle$  which  is not a staircase
relation when restricted to any sos $(a_0,\dots,a_{m-1})$ of $(e_i)_{i=0}^n$. Let $\mU$ be a
non-principal ultrafilter on $\N$, and define the  equivalence relation $\rel$ on $FIN_k$ by
$$s\rel
t \text{ if and only if } \conj{n}{s R_n t}\in \mU,$$
 where
$R_n=\rel_n\cup (\fin_k)^2$ is an equivalence relation on $\fin_k$. It is easy to see that $\rel$
is an equivalence relation. By  Theorem \ref{maint}, there is some sos $A=(a_n)_n$ on which $\rel$
is a staircase equivalence relation, say $\rel_{can}$. Choose $n$ large  enough such that:
\begin{enumerate}
\item $(a_0,\dots,a_{m-1})\pe (e_i)_{i=0}^n$
\item For $s,t\in \langle a_0,\dots,a_{m-1} \rangle$ one has that $s\rel t$ iff $s\rel_n t $.
\end{enumerate}
This can be done as follows: For every pair $s,t\in \langle a_0,\dots,a_{m-1} \rangle$, let
\begin{equation}
A_{s,t}=\left\{
\begin{array}{ll}
\conj{n}{s \rel_n t}\in \mU & \text{if } s \rel t \\
\conj{n}{s \nrel_n t}\in \mU & \text{if } s \nrel t.
\end{array} \right.
\end{equation}
Let $n=\min \bigcap_{s,t \in \langle a_0,\dots,a_{m-1} \rangle }A_{s,t} $. Then $\rel_n$ is $\rel$
restricted to $(a_0,\dots,a_{m-1})$, and hence is a staircase equivalence relation, a
contradiction.
\fprue
 \cor For every $m$
there is some $n=n(m)$ such that  for every equivalence relation $\rel$
 on $\langle e_0,\dots,e_n \rangle$ there is some sos $(a_0,\dots,a_{m-1})\pe( e_0,\dots,e_n)$
  such that $\rel$ is a canonical equivalence relation on $\langle a_0,\dots,a_{m-1} \rangle$.\qed
\fcor
\cor For every $m$ there is some $n=n(m)$ such that for every $(b_0,\dots,b_n)$ and  every
equivalence relation $\rel$
 on $\langle b_0,\dots,b_n \rangle$ there is some sos $(a_0,\dots,a_{m-1})\pe (b_0,\dots,b_n)$
  such that $\rel$ is a staircase equivalence relation when restricted to $\langle a_0,\dots,a_{m-1} \rangle$.
\fcor
\prue
Let $n=n(m)$ be given by Theorem \ref{finite}. Fix $b_0,\dots,b_n$, and an  equivalence relation
$\rel$. Let $F$ be the canonical isomorphism between $\langle e_0,\dots,e_n \rangle$ and $\langle
b_0,\dots,b_n \rangle$. Define $\rel'$ on $\langle e_i \rangle_{i=1}^n$ via $F$, i.e., $s\rel' t$
\iff $F(s)\rel F(t)$. Fix an sos $(c_i)_{i=0}^{m-1}\pe (e_i)_{i=1}^n$ and a staircase equivalence
relation $\rel_{can}$ such that $s\rel' t$ \iff $s \rel_{can} t$, for every $s,t\in \langle c_i
\rangle_{i=0}^{m-1}$. Let $b_i=F c_i$, for every $i=0,\dots,m-1$. Then $(b_0,\dots,b_m)$ is an sos
since sos are preserved under isomorphisms, and  $\rel_{can}$ is well defined on
$(b_0,\dots,b_{m-1})$. Since $\rel_{can}$ is staircase one has that
$$\text{ $s\rel_{can}t$ \iff
$F^{-1}s\rel_{can}F^{-1}t$}$$
 for every $s,t\in \langle
b_0,\dots,b_{m-1} \rangle$. Hence,
$$\text{$s\rel_{can}t$ iff $F^{-1}s\rel_{can} F^{-1}t$ iff $F^{-1}s\rel'
F^{-1}t$ iff $s\rel t$.}$$
 Therefore $\rel_{can}$ and $\rel$ coincide on $\langle
b_0,\dots,b_{m-1} \rangle$.
\fprue

\defi
We say that a staircase relation $\rel$ is \emph{symmetric} iff $I_1(\rel)=I_0(\rel)=I$,
$J_1(\rel)=J_0(\rel)=J$ and $l_j^{(0)}=l_j^{(1)}$ for every $j\in J$.
\fdefi

\cor\label{eojrghheruhgtr}
 For every $m$ there is some $n=n(m)$ such that  for every equivalence relation $\rel$
 on $\langle e_0,\dots,e_n \rangle$ there are \emph{disjointly supported} sos's $a_0,\dots,a_{m-1} \in \langle
(e_i)_{i=0}^n\rangle$
  such that $\rel$ is a symmetric staircase relation in $\langle a_0,\dots,a_{m-1} \rangle$.
\fcor
Before we give the proof of this, let us observe that if $a_0,\dots,a_{m-1}$ are  disjointly
supported $k$-vectors then the mapping $e_i\mapsto a_i$ extends to a lattice-isomorphism from
$\langle (e_i)_{i=0}^{m-1}\to \langle (a_i)_{i=0}^{m-1}\rangle$ that preserves the operation $T$.
\prue
Fix an integer $m$. Let $n$ be given by Theorem \ref{finite} when applied to $2m$. Suppose that
$\rel$ is an equivalence relation on $\langle (e_i)_{i=0}^n\rangle$. Then there is some sos
$(b_i)_{i=0}^{2m-1}$  such that $\rel$ is a staircase relation when restricted to $\langle
(b_i)_{i=0}^{2m-1}\rangle$.  Let $a_i=b_i+b_{2m-i-1}$ for every $0\le i\le m-1$. A typical vector
$b\in \langle (a_i)_{i=0}^{m-1}\rangle$ is of the form
$$b=\sum_{i=0}^{m-1}T^{k-r_i}b_i + \sum_{i=0}^{m-1}T^{k-r_i}b_{2m-i-1}.$$ Let $s,t\in
\langle(a_i)_{i=0}^{m-1}\rangle$. Then one has that
\begin{align}
\label{othriohhfuwww1}\min_i(s)=\min_i(t)& \text{ iff } \max_i(s)=\max_{i}(t),\text{ and}\\
\label{othriohhfuwww2}\theta_{i,l}^0(s)=\theta_{i,l}^0(t) &\text{ iff }
\theta_{i,l}^1(s)=\theta_{i,l}^1(t).
\end{align}
Let $(I_0,J_0,(l^{(0)}_j)_{j\in J_0},I_1,J_1,(l^{(1)}_j)_{j\in J_1},l^{(2)}_k)$ be the values of
$\rel$ when restricted to $\langle (b_i)_{i=0}^{2m-1}\rangle$. Using \eqref{othriohhfuwww1} and
\eqref{othriohhfuwww2} it follows that our fixed relation $\rel$ is when restricted to
$\langle(a_i)_{i=0}^{m-1}\rangle$ a symmetric staircase relation with values
$$(I_0\cup I_1,J_0\cup J_1,(l_j)_{j\in J_0\cup J_1},I_0\cup I_1,J_0\cup J_1,(l_j)_{j\in J_0\cup J_1},l^{(2)}_k)$$
and where for each $j\in J_0\cup J_1$
$$l_j=\left\{\begin{array}
{ll} \min \{l_j^{(0)},l_j^{(1)}\}& \text{if $j\in J_0\cap J_1$}\\
l_j^{(0)}&\text{if $j\in J_0\setminus J_1$}\\
l_j^{(1)}&\text{if $j\in J_1\setminus J_0$}.
\end{array}\right.$$
\fprue

\nota
\noindent (1) Pr\"{o}mel and Voigt were the firsts to observe  in \cite{pro-vo}  the Corollary
\ref{eojrghheruhgtr} for $\fin$. We thank the referee for pointing us out  this.

\noindent (2) Let $\mc S_k$ be the set of symmetric staircase relations of $\fin_k$, and set
$s_k=|\mc S_k|$. Using the notation from the Section \ref{counting} one has that
$$\mc S_k=\mc C_k \cup \conj{\rel \cap \rel_{\theta_l^{2}}}{\rel \in \mc A_k\setminus\mc C_k ,\text{ and $l=-1,1,\dots,k$}}.$$
Hence
$$s_k=c_k+(a_k-c_k)(k+1)=a_{k-1}+(k+1)(a_k-a_{k-1})=(k+1)!e_k(1)-k! e_{k-1}(1).$$
\fnota

\section{Canonical relations and continuous maps on $PS_{c_0}$}
Our result on equivalence relations on $\fin_k$ gives some consequences about equivalence relations
on $PS_{c_0}$. Let us start with some natural definitions.

For a fixed $\de>0$, let $k$ be the first integer such that $1/(1+\de)^{k-1}<\de$, and set
$\de_i=(1+\de)^{i-k}$, for $0\le i\le k$. For $0\le i \le k+1$, let
\begin{equation*}
\ga_i(\de)=\left\{\begin{array}{ll} \frac{\de_{i-1}+\de_i}2=\frac{\vep^{k-i}(\vep+1)}2 &\text{if } 1\le i \le
k \\
0 & \text{if }i=0 \\
\de_k=1 & \text{if }i=k+1
\end{array} \right.
\end{equation*}
and for $0\le i \le k$, let
\begin{equation*}
I_i^{(\de)}=\left\{\begin{array}{ll} \text{$[\ga_i(\de),\ga_{i+1}(\de))$} &\text{if } 0\le i <
k \\
\text{$[\ga_k(\de),\ga_{k+1}(\de)=1]$} & \text{if }i=k.
\end{array} \right.
\end{equation*}
We have then  that $\de_i\in I_i^{(\de)}$ for every $0\le i \le k$, and that
$[0,1]=\bigcup_{i=0}^kI_i^{(\de)}$, a disjoint union.

For $x=(x_m)_m\in PB_{c_0}$  and $n\in \N$, let $\Ga_n^{(\de)}(x)$ be the unique $0\le i \le k$
such that $x_n\in I_i^{(\de)}$, and define $\Ga_\de:PB_{c_0}\to \fin_{\le k}$ by
$\Ga_\de(x)=(\Ga^{(\de)}_n(x))_n$. Notice that $\Ga_\de (PS_{c_0})\con \fin_k$. A vector $x\in
PS_{c_0}$ is called a \emph{$\de$-sos} iff $\Ga_{\de}x $ is an sos. A block sequence $(x_n)_n$ of
vectors of $PS_{c_0}$ is called a \emph{$\de$-sos} iff every $x\in PS_{X}$ is a $\de$-sos.  The
next proposition is not difficult to prove.
\prop Fix $\rho\in [0,1]$, $x,y \in PB_{c_0}$, and a
$k$-vector $s$ of $\fin_k$.  Let $i$  be the unique integer such that $\rho\in I_i^{(\de)}$. Then,
\begin{enumerate}
\item $\Ga_\de (x+y)=\Ga_\de(x)+\Ga_\de (y)$ and $\Ga_\de (\rho e_n)=T^{k-i}\Ga_\de
(e_n)=T^{k-i}(\Theta_\de^{-1}e_n ) $.
\item $\Ga_\de (\rho \Theta_\de^{-1}x)=T^{k-i}\Ga_\de
(\Theta_\de^{-1}x)$. It follows that  if $(a_n)_n$ is an sos $k$-block sequence, then
$(\Theta_\de^{-1}a_n)_n$ is a $\de$-sos.\qed
\end{enumerate}

\fprop

\defi
Given a staircase mapping $f$ of $\fin_k$, we consider the following two extensions to an arbitrary $\de$-sos
$X=(x_n)_n$. The fist one is $f^{(0)}:PS_X\to \fin_{\le k}$, closing the following diagram:

\begin{center}
\begin{figure}[h]
\includegraphics[scale=1]{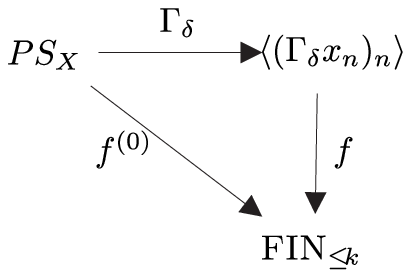}
\end{figure}
\end{center}
The second one is $f^{(1)}:PS_X\to PB_{c_0}$, defined by $f^{(1)}(x)(n)=x(n)$ iff $f^{(0)}x(n)\neq
0$.
\fdefi
\prop Fix a staircase $f$, and some $\de$-sos $X$.
\begin{enumerate}
\item $(f\odot g)^{(i)}=f^{(i)}\odot g^{(i)}$, for $i=0,1$ and $\odot$ equal to $\vee$ or $\wedge$.
\item  $f^{(1)}$ is a Baire class 1 function.
\item If $f^{(1)}x=f^{(1)}y$, then
$f^{(0)}x=f^{(0)}y$ for every $x,y\in PS_X$.
\item $\|\Theta_\de^{-1}f^{(0)}x-f^{(1)}x\|\le \de$ for every $x\in PS_X$.
\item
For every $k$-vector $a\in \langle (\Ga_\de x_n)_n\rangle$, $f^{(1)}\Theta^{-1}_\de a=f^{(0)}\Theta^{-1}_\de
a=fa$. Therefore $f^{(1)}\Theta^{-1}_\de a=f^{(1)}\Theta^{-1}_\de b$ iff $f^{(0)}\Theta^{-1}_\de
a=f^{(0)}\Theta^{-1}_\de b$, for every $k$-vectors $a,b\in \langle (\Ga_\de x_n)_n\rangle$.
\item For every $x\in PS_X$ there is some $k$-vector $\bar{x}$ such that
$\|x-f^{(0)}\Theta^{-1}_\de \bar{x}\|\le \de$ and  $f^{(0)}x=f^{(0)}\Theta^{-1}_\de
\bar{x}$.

\end{enumerate}
\fprop
\prue
(i) is not difficult to check. Let us show (ii). To do this, suppose that $f$ is a staircase
mapping. Then $f$ is in the algebraic closure of $\mathcal{F}$ (see Definition \ref{defoff}), i.e.,
there is a finite list $f_0,\dots,f_n\in \mathcal{F}$ such that $f= f_0 \odot_0 f_1 \odot_1 f_2
\odot_2 \cdots \odot_{n-1} f_n$, where $\odot_i$ is either $\vee$ or $\wedge$ for every
$i=0,\dots,n-1$. By point (i) one has that $f^{(1)}= f_0^{(1)} \odot_0 f_1^{(1)} \odot_1 f_2^{(1)}
\odot_2 \cdots \odot_{n-1} f_n^{(1)}$. Since for every point $x\in PS_X$ the support of
$f_i^{(1)}(x)$ is finite,  we may assume that $f\in \mathcal{F}$. We give the proof for the case
$f=\min_i$. The other cases  can be shown in a similar way. For $l>0$ we define the following
perturbations of the intervals $I^{(\de)}_i$, let
\begin{equation*}
I_{i,l}^{(\de)}=\left\{ \begin{array}{ll} (\ga_i(\de)-\frac1l,\ga_{i+1}(\de)) & \text{if } i<k \\
(\ga_k(\de)-\frac1l,1] &\text{if }i=k.
\end{array}\right.
\end{equation*}
These are open intervals of $PS_{c_0}$. For each $l$, let $f_l:PS_X\to PB_X$ be defined for $n\in \N$ as
follows,
\begin{equation*}
f_l(x)(n)=\left\{ \begin{array}{ll} x(n)  & \text{if $x(n)\in I_{i,l}^{(\de)}$ and for all $m<n$ $x(n)\in [0,\ga_i(\de))$}  \\
0 &\text{if not.}
\end{array}\right.
\end{equation*}
Let us see that $f_l$ is continuous, and that $f_l\to_l f$. Suppose that $x_r\to_r x$, with $x_r, x\in PS_X$.
Let $n$ be the unique integer such that $f_l(x)(n)=x(n)>0$,  i.e., $x(n)\in I_{i,l}^{(\de)}$  and $x(m)\in
[0,\ga_i(\de))$ for every $m<n$. Since both sets are open, there must be some $r'$ such that $x_{r''}(n)\in
I_{i,l}^{(\de)}$ and $x_{r''}(m)\in [0,\ga_i(\de))$, for every $r''>r'$  and every $m<n$. Therefore, for all
$r''>r'$, $f_l(x_{r''})=f_lx$.  Let us check now that $f_l\to f$. Fix $x$, and we  show that $f_l(x)\to
f(x)$. Again, Let $n$ be the unique integer such that $f_l(x)(n)=x(n)>0$. Let $l'$ be such that $x(m)\in [0,
\ga_i(\de)-1/l)$ for every $m<n$. Then  $f_{l''}x (m)=0$ and $f_{l''}x(n)=x(n)$,  for every $l''\ge l'$ and
every $m<n$. Also,  $f_{l''}x(m)=0$ for every $m>n$. All this implies that $f_{l''}(x)=f(x)$.

The rest of the points (iii)-(vi) are not difficult to prove. We leave the details to the reader.
 \fprue
 For
an equivalence relation $R$, and $x\in PS_{c_0}$, the \emph{$R$-equivalence class} of $x$ is denoted by
$[x]_R$.
\prop Fix $\de>0$, a staircase equivalence relation $R_f$, and  a $k$-block sequence
$A=(a_n)_n$, where $k=k(\de)$. Set $X=(x_n=\Theta^{-1}_\de a_n)_n$ and $R=R_{f^{(1)}}$.
\begin{enumerate}
\item For every $x\in PS_X$ there is a $k$-vector $\bar{x}$ of $A$ such that $\|x-\Theta^{-1}_\de
\bar{x}\|\le \de$ and $[x]_{R}\con ([\Theta^{-1}_\de\bar{x}]_R)_\de$.
\item For every $x,y\in PS_X$, if $(x,y)\in R$, then $(x,z)\in R$,  for every $x\wedge y \le_L z \le_L x \vee y
$.
\end{enumerate}

\fprop
\prue
To prove (i), fix  $x\in PS_X$, and let $\bar{x}$ be a $k$-vector of $A$ such that $\|x-\Theta^{-1}_\de
\bar{x}\|\le \de$ and $f^{(0)}x=f^{(0)}\Theta^{-1}_\de \bar{x}$. Set  $x'=\Theta^{-1}_\de \bar{x}$. We  show
that $[x]_R\con ([x']_R)_\de$. Suppose that $y\in [x]_R\cap PS_X$. Then $f^{(1)}x=f^{(1)}y$, and hence
$f^{(0)}y=f^{(0)}x=f^{(0)}x'$. Let $\bar{y}$ be a $k$-vector of $A$ such that $\|y-\Theta^{-1}_\de
\bar{y}\|\le \de $
 and $f^{(0)}y=f^{(0)}\Theta^{-1}_\de \bar{y}$, and set $y'=\Theta^{-1}_\de \bar{y}$.
Then, $f^{(0)}x'=f^{(0)}y'$, which implies that $f^{(1)}x'=f^{(1)}y'$, i.e., $y'\in [x']_R$ and
hence $y\in ([x']_R)_\de$.

(ii): By Proposition \ref{stairr}, we may assume that $f\in \mathcal{F}$. Again, we give a proof for the case
$f=\min_i$, since the other cases can be shown in a similar way. Suppose that $(x,y)\in R_{f^{(1)}}$, and fix
$z\in PS_X $ with $x\wedge y \le_L z \le_L x \vee y $. Let $n$ be the unique integer such that
$f^{(1)}x(n)=x(n)=y(n)=f^{(1)}y(n)>0$.  Then $x(m),y(m)\in [0,\ga_i(\de))$ for every $m<n$. Therefore,
$z(n)=x(n)=y(n)$  and $z(m)\in [0,\ga_i(\de))$ for every $m<n$. This implies that $f^{(1)}(z)=f^{(1)}(x)$.
\fprue
\defi
A \emph{$\de$-staircase equivalence relation} is $R_{f^{(1)}}$ for some staircase $f$.
\fdefi

The next result is the interpretation of  Theorem \ref{maint} in terms of equivalence relations on
$PS_X$.

\prop\label{gen88}
Let $R$ be an equivalence relation on $PS_X$. Then for every $\de>0$ there is some $\de$-sos $X$ and some
$\de$-staircase equivalence relation $\widetilde{R}$ such that:
\begin{enumerate}
\item $R$ and $\widetilde{R}$ coincide in an $\vep$-net of $PS_X$ for some $\vep<\de$.
\item  For every $\widetilde{R}$-class $\al$ on
$PS_X$ there is a $R$-class $\be$ on $PS_X$ such that $\al \con \be_\de$.
\end{enumerate}

\fprop
\prue

Fix $\de$, and let $k=k(\de)$. Define $\bar{R}$ on $\fin_k$ via $\Theta_\de$.  Then there is some sos
$k$-block sequence $A=(a_n)_n$ and some staircase equivalence relation $R_f$ such that $\bar{R}$ and $R_f$
coincide on $\langle A \rangle$. Set $\widetilde{R}=R_{f^{(1)}}$ and $X=(x_n)_n$, where
$x_n=\Theta_\de^{-1}a_n$ for every $n$.

\noindent (i): For $\vep=(1+\de)^{k-1}$, $N=\Theta_\de^{-1}(\langle( a_n)_n\rangle)$ is a
$\vep$-net of $PS(X)$ satisfying our requirements.

\noindent (ii): For a fixed  $x\in PS_X$ choose some $k$-vector $\bar{x}$ of $A$ such that
$\|x-x'\|\le \de$ and $f^{(0)}x=f^{(0)}x'$, where  $x'=\Theta^{-1}_\de \bar{x}$. We  show that
$[x]_{\widetilde R}\con ([x']_R)_\de$. Suppose that $y\in PS_X$ is such that $f^{(1)}x=f^{(1)}y$.
Pick some $k$-vector $\bar{y}$ of $A$  such that  $\|y-y'\|\le \de$ and $f^{(0)}y=f^{(0)}y'$ where
$y'=\Theta^{-1}_\de \bar{y}$. Then, $f^{(0)}x=f^{(0)}y$ and hence $f^{(0)}x'=f^{(0)}y'$, which
implies that $f^{(1)}x'=f^{(1)}y'$. Therefore, $y'\in [x']_R$.
\fprue

In the case of equivalence relations with some additional properties, we have the following stronger result.
\prop
Fix $\de, \ga>0$, set $k=k(\de)$, and suppose that $R$ is an equivalence relation on $PS_{c_0}$ such that
\begin{enumerate}
\item for every $x,y\in PS_{c_0}$ and every $z\in PS_{c_0} $ with $x\wedge y \le_L z \le_L x \vee y $, if
$(x,y)\in R$, then   $(x,z)\in R$, and
\item  for every sos $k$-block sequence $B=(b_n)_n$  and every $x\in
PS_{(\Theta^{-1}_{\de}b_n)_n}$ there is some $k$-vector $\bar{x}$ of $B$ such that
 $[x]_R\con
([\Theta^{-1}_{\de}\bar{x}]_R)_\ga$.
\end{enumerate}
Then, there is some $\de$-sos $X$  and some $\de$-staircase equivalence relation $\widetilde{R}$ such that

\noindent \emph{(a)} for every $R$-equivalent classes $\al$ in $PS_X$, there is a $\widetilde{R}$-equivalent
class $\be$ in $PS_X$ such that $\al\con \be_{\de+\ga}$, and

\noindent \emph{(b)} for every $\widetilde{R}$-equivalence class $\be$ there is a $R$-equivalence class $\al$
such that $\be\con (\al)_\de$.
\fprop
\prue

Define $\bar{R}$ on $\fin_k$ via $\Theta_{\de}$. Then, there is some sos $A=(a_n)$ and some
staircase equivalence relation $R_f$ such that $\bar{R}$ is $R_f$ on $\langle  A \rangle$. Let
$\widetilde{R}=R_{f^{(1)}}$, and  $X=(x_n)_n$, where $x_n=\Theta_{\de}^{-1}a_n$ for every $n$.
\emph{(b)} is shown in Proposition \ref{gen88}. Let us show \emph{(a)}.  Fix $x\in PS_X$, and
choose a $k$-vector $\bar{x}$ of $A$ such that $[x]_R\con ([x']_R)_\ga$ where $x'=\Theta^{-1}_\de
\bar{x}$. Let us show that $[x']_R\con ([x']_{\widetilde{R}})_\de$ on $PS_X$. Fix $y\in [x']_R$.
Then, there is some $k$-vector $\bar{y}$ of $A$ such that $x'\wedge y \le_L y'\le_L x' \vee y$ and
$\|y-y'\|\le \de$, where $y'=\Theta^{-1}_\de \bar{y}$. Hence, $y'\in [x']_R$, and therefore,
$y'\in [x']_{R'}$.
\fprue

\noindent \textbf{Acknowledgement } I wish to thank Stevo Todorcevic for giving me the problem, as
well as for continuous support during the course of this work.

%\vspace{1cm}
%
%
%\begin{tabular}{l}
%Jordi L\'{o}pez-Abad \\
%Universit\'{e} Paris 7- Denis Diderot\\
%C.N.R.S. - UMR 7056 \\
%2, Place Jussieu- Case 7012\\
%75251 Paris Cedex 05 \\
%France \\
%{\em E-mail}: \texttt{abad@logique.jussieu.fr}
%\end{tabular}

\end{document}